\documentclass[12pt]{article}
\usepackage[dvips]{graphicx}
\usepackage{amsmath,amsthm,amsfonts,amssymb,amscd,epic,eepic}

\newtheorem{thm}{Theorem}[section]
\newtheorem{cor}[thm]{Corollary}
\newtheorem{defn}[thm]{Definition}

\newtheorem{lem}[thm]{Lemma}
\newtheorem{prop}[thm]{Proposition}

\newtheorem{rem}[thm]{Remark}
\newtheorem{rems}[thm]{Remarks}

\title{$(2+1)$-dimensional topological quantum field theory 
from subfactors and Dehn surgery formula for 3-manifold invariants}

\author{
{\sc Yasuyuki Kawahigashi}\footnote{Supported in part by the
Grants-in-Aid for Scientific Research, JSPS.}\\
Department of Mathematical Sciences\\
University of Tokyo, Komaba, Tokyo, 153-8914, JAPAN\\
e-mail: {\tt yasuyuki@ms.u-tokyo.ac.jp}\\
\vphantom{X}\\
{\sc Nobuya Sato}$^*$ \\
Department of Mathematics and Information Sciences\\
Osaka Prefecture University, Sakai, Osaka, 599-8531, JAPAN\\
e-mail: {\tt nobuya@mi.cias.osakafu-u.ac.jp}\\
\vphantom{X}\\
{\sc Michihisa Wakui}\\
Department of Mathematics\\
Osaka University, Toyonaka, Osaka, 560-0043, JAPAN \\
e-mail: {\tt wakui@math.sci.osaka-u.ac.jp}}

\begin{document}

\maketitle

\begin{abstract}
In this paper, we establish the general theory of $(2+1)$-dimensional 
topological quantum field theory (in short, TQFT) with a Verlinde basis. 
It is a consequence that we 
have a Dehn surgery formula for 3-manifold invariants for this kind of 
TQFT's. We will show that Turaev-Viro-Ocneanu unitary TQFT's obtained from 
subfactors satisfy the axioms of TQFT's with Verlinde bases. Hence, in a 
Turaev-Viro-Ocneanu TQFT, we have a Dehn surgery formula for 3-manifolds. 
It turns out that this Dehn surgery formula is nothing but the formula of 
the Reshetikhin-Turaev invariant constructed from a tube system, 
which is a modular category corresponding to the quantum double construction 
of a $C^*$-tensor category. In the forthcoming paper, we will exbit 
computations of Turaev-Viro-Ocneanu invariants for several \lq\lq basic 
3-manifolds ''. In Appendix, we discuss the relationship 
between the system of $M_\infty$-$M_\infty$ bimodules arising from the 
asymptotic inclusion $M \vee M^{op} \subset M_\infty$ constructed from 
$N \subset M$ and the tube system obtained from a subfactor $N \subset M$.
\end{abstract}

\section{Introduction}
Since the discovery of the celebrated Jones polynomial \cite{J}, 
there have been a great deal of studies on quantum invariants of knots, 
links, and 3-manifolds.  Among them, we deal with the two constructions
in this paper; Ocneanu's generalization of the Turaev-Viro 
invariants \cite{TuraevViro} of 3-manifolds based on triangulation
and the Reshetikhin-Turaev invariants \cite{RT} of those based
on Dehn surgery.  The former uses quantum $6j$-symbols arising from
subfactors and the tensor categories we have are not necessarily braided,
while the the latter invariants require braiding, or more precisely,
modular tensor categories.  Methods to construct a modular tensor
category from a rational tensor category have been studied by
several authors and are often called ``quantum double'' constructions.
In the setting of subfactor theory, such a method was first found
by Ocneanu using his generalization of the Turaev-Viro topological
quantum field theory (TQFT) and Izumi \cite{Izumi1} later gave a formulation
based on sector theory.

We start with a (rational unitary) tensor category (arising from
a subfactor), which is not braided in general.  
Then we obtain two TQFT's with the methods above as follows.  One is
a Turaev-Viro-Ocneanu TQFT based on a state sum and triangulation which
directly uses $6j$-symbols of the tensor category.  The other is
a Reshetikhin-Turaev TQFT arising from the modular tensor category
we obtain with the quantum double construction of the original category.
It is natural to ask what relation we have between these two
TQFT's.  In order to make a general study on such TQFT's, we first
establish a Dehn surgery formula (Proposition \ref{Proposition1}) 
for a general TQFT, which represents an invariant of
a 3-manifold in terms of a weight sum of invariants of links that
we use for the Dehn surgery construction of the 3-manifold.  This
formula uses a basis of the Hilbert space for $S^1\times S^1$
arising from the TQFT.  We show that several nice properties for
such a basis are mutually equivalent (Theorem \ref{Theorem6}) and we say 
that a basis is a {\it Verlinde basis} when such properties hold 
(Definition \ref{vb}).  We next show that 
a Turaev-Viro-Ocneanu TQFT has a Verlinde basis in this sense 
(Theorem \ref{verlinde-basis}) and then as a corollary of the Dehn surgery 
formula, we conclude that the above two TQFT's are identical for 
any closed 3-manifold (Theorem \ref{main}). 
(Actually, the claim that a Turaev-Viro-Ocneanu TQFT has a Verlinde
basis has been announced by Ocneanu and a proof is presented in
\cite[Chapter 12]{EK}, but normalizations are inaccurate
there, so we include a proof for this claim for the sake of completeness
here, along the line of Ocneanu's tube algebra and braiding arising 
from the Turaev-Viro-Ocneanu TQFT.)
This identity result has been announced by two of us in \cite{SatoWakui}.
It has been also announced in page 244 of Ocneanu \cite{O},
but it seems to us that his line of
arguments are different from ours which relies on our
Dehn surgery formula.  This result also proves a conjecture
in \cite[Section 8.2]{Mueger} and gives an answer to a question
in \cite[page 546]{Kerler}. 
We refer to the book \cite{EK} for subfactor theory and its applications.
Throughout this paper, we only consider subfactors with 
finite depth and finite index. 

\bigskip
\noindent
Acknowledgement: A part of this work was done while the first two 
authors visited Mathematical Sciences Research Institute at Berkeley. 
We appreciate their financial supports and hospitality.

\section{$(2+1)$-dimensional TQFT with Verlinde basis and Dehn surgery 
formula}

Giving a $(2+1)$-dimensional topological quantum field theory $Z$, we have 
a 3-manifold invariant $Z(M)$ in the canonical way. The purpose of this 
section is to describe $Z(M)$ in terms of the $S$-matrix and framed links, 
and to give a criterion for the Verlinde identity \cite{Verlinde} to hold in 
our framework. 
\par 
Throughout this section, we use the following notations. 
For an oriented manifold $M$, we denote by $-M$ the same manifold with the opposite orientation. For an orientation preserving diffeomorphism $f:M\longrightarrow N$, we denote by $-f$ the same diffeomorphism viewed as an orientation preserving diffeomorphism $-M\longrightarrow -N$. We regard the empty set $\emptyset $ as an oriented closed surface with the unique orientation. The dual vector space of a finite dimensional vector space $V$ over $\mathbb{C}$ is denoted by $V^{\ast }$, and the dual basis of a basis $\{ v_i\} _i$ is denoted by $\{ v_i^{\ast }\} _i$. 
The closure of a subspace $X$ in the $3$-sphere $S^3$ is denoted by $\overline{X}$. 

\par \medskip \noindent 
{\bf Convention of orientations for manifolds}.\ 
\par
\begin{figure}[hbtp]

\begin{center}
\includegraphics[width=5cm]{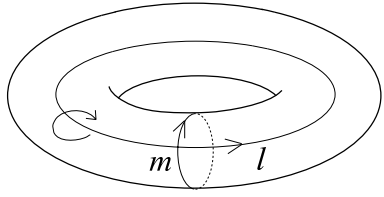} 
\caption{the orientation of torus $S^1\times S^1$ \label{Figure1}}
\end{center}
\end{figure}

 
Throughout this section, we assume that the $3$-sphere $S^3$ is oriented, and assume that any $3$-dimensional submanifold of $S^3$ is oriented by the orientation induced from $S^3$. We also assume that the solid torus $D^2\times S^1$ equips with the orientation such that the diffeomorphism $h:D^2\times S^1\longrightarrow \mathbb{R}^3\subset S^3$ defined by
$$h((x,y),e^{i\theta })=((2+x)\cos \theta , (2+x)\sin \theta ,y)$$
is orientation preserving, and assume that the torus $S^1\times S^1$ is oriented by the orientation induced from $D^2\times S^1$ (see Figure \ref{Figure1}). 
 
\par \medskip 
Let us recall the definition of $(2+1)$-dimensional topological quantum field theory due to Atiyah \cite{Atiyah}, \cite{FSS}. 
By a $3$-cobordism we mean a triple $(M; \Sigma _1,\Sigma _2)$ consisting of a compact oriented $3$-manifold $M$ and two closed oriented surfaces $\Sigma _1, \Sigma _2$ such that $\partial M=(-\Sigma _1)\cup \Sigma _2$ and $\Sigma _1\cap \Sigma _2=\emptyset $. For an oriented closed surface $\Sigma $, the identity cobordism is given by $Id_{\Sigma }=(\Sigma \times [0,1];\Sigma \times \{ 0\} ,\Sigma \times \{ 1\} )$.

\par 
\bigskip 
\begin{defn}
Let $Z$ be a functor consisting of the following three functions. 
\begin{enumerate}
\item[(1)] 
To each closed oriented surface $\Sigma $, it assigns a finite dimensional $\mathbb{C}$-vector space $Z(\Sigma )$.
\item[(2)] 
To each cobordism $W$, it assigns a $\mathbb{C}$-linear map $Z_W$.
\item[(3)] 
To each orientation preserving diffeomorphism $f$ between closed oriented surfaces, it assigns a $\mathbb{C}$-linear isomorphism $Z(f)$. 
\end{enumerate}
If $Z$ has the following properties, it is called a {\rm $(2+1)$-dimensional 
topological quantum field theory} (TQFT in short).
\begin{enumerate}
\item[(i)] $Z$ is functorial with respect to the composition of orientation preserving diffeomorphisms. More precisely, 
\begin{enumerate} 
\item[(a)] $Z(g\circ f)=Z(g)\circ Z(f)$  
for orientation preserving diffeomorphisms $f$ and $g$. 
\item[(b)] $Z(id_{\Sigma })=id_{Z(\Sigma )}$ for an oriented closed surface $\Sigma $. 
\end{enumerate}
\item[(ii)] $Z$ is fuctorial with respect to the composition of $3$-cobordisms. More precisely, 
\begin{enumerate} 
\item[(a)] If two cobordisms $W_1=(M_1; \Sigma _1,\Sigma _2)$ and $W_2=(M_2; \Sigma _2,\Sigma _3)$ are obtained by cutting a cobordism $W=(M;\Sigma _1,\Sigma _3)$ along $\Sigma _2$ in $M=M_1\cup M_2$, then 
$$Z_W=Z_{W_2}\circ Z_{W_1}.$$
\item[(b)] $Z(i_1)^{-1}\circ Z_{Id_{\Sigma }}\circ Z(i_0)=id_{Z(\Sigma )}$ for the identity cobordism on $\Sigma $, where $i_t:\Sigma \longrightarrow \Sigma \times \{ t\} $ \ $(t=0,1)$ are orientation preserving diffeomorphisms defined by $i_t(x)=(x,t),\ x\in \Sigma $. 
\end{enumerate}
\item[(iii)] Let $W=(M,\Sigma _1,\Sigma _2),\ W'=(M',\Sigma _1',\Sigma 
_2')$ be two cobordisms. Suppose that there is an orientation preserving 
diffeomorphism $h:M\longrightarrow M'$ such that $h(\Sigma _i)=\Sigma 
_i'$ for $i=1,2$. Then, for $f_1:=-h\vert _{\Sigma _1}:\Sigma 
_1\longrightarrow \Sigma _1' \ and \ f_2:=h\vert _{\Sigma _2}:\Sigma _2\longrightarrow \Sigma _2'$, the following diagram commutes. 
$$\begin{CD}
Z(\Sigma _1)@>Z(f_1)>> Z(\Sigma _1') \\ 
@VZ_{W}VV @VVZ_{W'}V \\ 
Z(\Sigma _2) @>>Z(f_2)> Z(\Sigma _2')
\notag 
\end{CD}$$
\item[(iv)] Let $W_1=(M;\Sigma _1,\Sigma _2)$, $W_2=(N;\Sigma _2',\Sigma _3)$ be two cobordisms, and $f:\Sigma _2\longrightarrow \Sigma _2'$ an orientation preserving diffeomorphism. Then, for the cobordism $W:=(N\cup _fM; \Sigma _1, \Sigma _3)$ obtained by gluing of $W_1$ to $W_2$ along $f$, 
$$Z_W=Z_{W_2}\circ Z(f)\circ Z_{W_1}.$$
\item[(v)] There are natural isomorphisms 
\begin{enumerate}
\item[(a)] $Z(\emptyset )\cong \mathbb{C}$. 
\item[(b)] $Z(\Sigma _1\coprod \Sigma _2)\cong Z(\Sigma _1)\otimes Z(\Sigma _2)$ for oriented closed surfaces $\Sigma _1$ and $\Sigma _2$.
\item[(c)] $Z(-\Sigma )\cong Z(\Sigma )^{\ast }$ for an oriented closed surface $\Sigma $. 
\end{enumerate}
\end{enumerate}
\end{defn}

\par \medskip 
\begin{rems}
1. For an oriented closed $3$-manifold $M$, we have a cobordism $W=(M;\emptyset ,\emptyset )$. This cobordism $W$ induces a linear map 
$$\mathbb{C}\cong Z(\emptyset )\xrightarrow {Z_W}Z(\emptyset )\cong \mathbb{C}.$$
We denote by $Z(M)$ the image of $1$ under the above map. By the condition (iii) we see that $Z(M)$ is a topological invariant of $M$. 
\par 
\noindent 
2. Let $\varGamma_{\Sigma }$ denote the mapping class group of the oriented closed surface $\Sigma $. Then, we have a representation of $\varGamma_{\Sigma }$
$$\rho :\varGamma_{\Sigma }\longrightarrow GL(Z(\Sigma )),\ [f] \longmapsto Z(f),$$
where $[f]$ denotes the isotopy class of $f$. 
\par 
\noindent 
3. For an oriented closed surface $\Sigma $, we have $\dim Z(\Sigma )=Z(\Sigma \times S^1)$. 
\par 
\noindent 
4. By a $3$-cobordism with parametrized boundary we mean a triple $(M; j_1,j_2)$ consisting of a compact oriented $3$-manifold $M$, an orientation reversing embedding $j_1:\Sigma _1\longrightarrow \partial M$ and an orientation preserving embedding $j_2:\Sigma _2\longrightarrow \partial M$ such that $\partial M=(-j_1(\Sigma _1))\cup j_2(\Sigma _2)$ and $j_1(\Sigma _1)\cap j_2(\Sigma _2)=\emptyset $. Any cobordism with parametrized boundary ${\cal W}=(M; \Sigma _1\xrightarrow {\ j_1\ } \partial M,\Sigma _2 \xrightarrow {\ j_2\ } \partial M)$ induces a linear map $Z_{\cal W}:Z(\Sigma _1)\longrightarrow Z(\Sigma _2)$ such that the following diagram commutes.

\begin{equation*}
\begin{CD}
Z(\Sigma _1) @>Z_{\cal W}>> Z(\Sigma _2) \\ 
@VZ(j_1)VV @VVZ(j_2)V \\ 
Z(j_1(\Sigma _1)) @>Z_W>> Z(j_2(\Sigma _2)) 
\end{CD} 
\end{equation*}

Here, we set $W=(M;j_1(\Sigma _1),j_2(\Sigma _2))$. 
\end{rems}

\par \bigskip 
\par
\begin{figure}[hbtp]
\begin{center}
\scalebox{1.3}{\includegraphics[height=2.5cm]{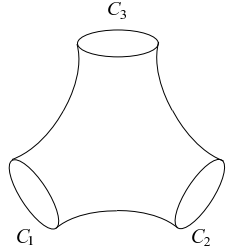} }
\caption{the compact oriented surface $Y$ \label{Figure5}}
\end{center}
\end{figure}

 \bigskip 
Let $Z$ be a $(2+1)$-dimensional TQFT. We
consider the cobordism $W:=(Y\times S^1; \Sigma _1\sqcup \Sigma _2, \Sigma _3)$, where $Y$ is
the compact oriented surface in $\mathbb{R}^3$
depicted in Figure \ref{Figure5} and $\Sigma _i=C_i\times S^1$ for $i=1,2,3$. 
Then, $W$ induces a 
linear map $Z_W: Z(S^1\times S^1)\otimes Z(S^1\times S^1)\longrightarrow
Z(S^1\times S^1)$. It 
can be easily verified that the map $Z_W$ gives an associative algebra structure
on $Z(S^1\times S^1)$. The identity element of this algebra is given by
$Z_{W_0}(1)$, where $W_0:=(D^2\times S^1; \emptyset , S^1\times S^1)$. We
call this algebra the {\it fusion algebra associated with $Z$}.

Let us introduce a Dehn surgery formula of $Z(M)$. 
Let $Z$ be a $(2+1)$-dimensional TQFT, and $S:S^1\times S^1\longrightarrow S^1\times S^1$ the orientation preserving diffeomorphism defined by $S(z,w)=(\bar{w},z),\ (z,w)\in S^1\times S^1$, where we regard $S^1$ as the set of complex numbers of absolute value $1$. Given a basis $\{ v_i\} _{i=0}^m$ of $Z(S^1\times S^1)$, we define $S_{ji}\in \mathbb{C}$ by $Z(S)v_i=\sum_{j=0}^mS_{ji}v_j$. 
\par 
Let $\{ v_i\} _{i=0}^m$ be a basis of $Z(S^1\times S^1)$. Then, we can define a framed link invariant as follows. Let $L=L_1\cup \cdots \cup L_r$ be a framed link with $r$-components in the 3-sphere $S^3$, and $h_i: D^2\times S^1\longrightarrow N(L_i)$ be the framing of $L_i$ for each $i\in \{ 1, \cdots , r\} $, where $N(L_i)$ denotes the tubular neighborhood of $L_i$. 
We fix an orientation for $\partial N(L_i)$ such that $j_i:=h_i\vert _{\partial D^2\times S^1}:S^1\times S^1\longrightarrow \partial N(L_i)$ is orientation preserving. Since the orientation for $N(L_i)$ is not compatible with the orientation for the link exterior $X:=\overline{S^3-N(L_1)\cup \cdots \cup N(L_r)}$, we can consider the cobordism with parametrized boundary ${\cal W}_L:=(X; \coprod \limits_{i=1}^r j_i ,\emptyset )$. 
This cobordism induces a $\mathbb{C}$-linear map $Z_{{\cal W}_L}: \bigotimes \limits_{i=1}^rZ(S^1\times S^1)\longrightarrow \mathbb{C}$. 
It is easy to see that for each $i_1,\cdots ,i_r=0,1,\cdots ,m$ the complex number $J(L; i_1, \cdots , i_r):=Z_{{\cal W}_L}(v_{i_1}\otimes \cdots \otimes v_{i_r})$  
is a framed link invariant of $L$. 

\par\bigskip 
\begin{prop}
\label{Proposition1}
Let $Z$ be a (2+1)-dimensional TQFT and  $\{ v_i\} _{i=0}^m$ a basis of $Z(S^1\times S^1)$ such that $v_0$ is the identity element in the fusion algebra. 
Let $M$ be a closed oriented 3-manifold obtained from $S^3$ by Dehn surgery along a framed link $L=L_1\cup \cdots \cup L_r$. 
Then, the 3-manifold invariant $Z(M)$ is given by the formula
$$Z(M)=\sum_{i_1, \cdots , i_r=0}^m S_{i_1,0}\cdots S_{i_r,0}J(L;i_1, \cdots , i_r).$$ 
\end{prop}

\medskip \noindent 
{Proof}.\ 
Let $X$ be the link exterior of $L$ in $S^3$, and $h_i:D^2\times S^1\longrightarrow \partial N(L_i)$ the framing of $L_i$ for each $i\in \{ 1,\cdots ,r\} $. 
Then, by using attaching maps $f_i: S^1\times S^1\longrightarrow \partial N(L_i)$ satisfying $f_i(S^1\times 1)=h_i(1\times S^1)$, $i=1,\cdots ,r$, we have 
$$M=X\bigcup \nolimits_{\amalg _{i=1}^rf_i}(\coprod^rD^2\times S^1).$$
Therefore, 
$$Z(M)=Z_{W_2}\circ Z(\coprod _{i=1}^rf_i)\circ Z_{W_1},$$
where $W_1:=(\coprod^r(D^2\times S^1); \emptyset , \coprod ^r(S^1\times S^1))$ and $W_2:=(X; -\partial X, \emptyset )$. This implies that 
$$Z(M)=Z_{\mathcal{W}_L}\circ (\bigotimes _{i=1}^rZ(j_i^{-1}))\circ (\bigotimes _{i=1}^rZ(f_i))\circ (\bigotimes ^rZ_{W_0}),$$
where $j_i:=h_i\vert _{\partial D^2\times S^1}:S^1\times S^1\longrightarrow \partial N(L_i)$ for $i=1,\cdots ,r$. 
\par 
For each $i$ we can choose $f_i$ satisfying $j_i^{-1}\circ f_i=S$ up to isotopy 
as for $f_i$ satisfying $f_i(S^1\times 1)=h_i(1\times S^1)$.    
Then, we have 
$$Z(M)=Z_{\mathcal{W}_L}\circ (\bigotimes _{i=1}^rZ(S))\circ (\bigotimes ^rZ_{W_0}).$$
This implies that 
$$Z(M)=\sum_{i_1, \cdots , i_r=0}^m S_{i_1,0}\cdots S_{i_r,0}J(L;i_1, \cdots , i_r).$$ 
This completes the proof. \hfill Q.E.D. 
\par \bigskip 

\begin{figure}[hbtp]
\begin{center}
\includegraphics[width=4cm]{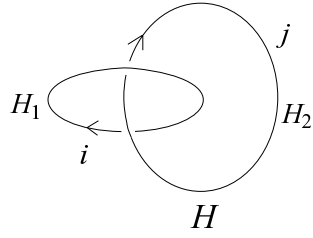} 
\caption{the Hopf link $H$ \label{Figure2}}
\end{center}
\end{figure}

\par \bigskip 
\begin{prop}
\label{Proposition2}  
Let $Z$ be a $(2+1)$-dimensional TQFT and $\{ v_i\} _{i=0}^m$ a basis of $Z(S^1\times S^1)$. 
Let $H$ be the Hopf link depicted as in Figure \ref{Figure2}, and $U$ the orientation preserving diffeomorphism from $S^1\times S^1$ to $-S^1\times S^1$ defined by $U(z,w)=(z,\bar{w})$. 
Then, the following are equivalent. 
\begin{enumerate}
\item[(i)] $J(H;i,j)=S_{ij}$ for all $i,j=0,1,\cdots ,m$. 
\item[(ii)] $Z(U)v_i=v_{i}^{\ast }$ for all $i=0,1,\cdots ,m$.
\end{enumerate} 
\end{prop}

\begin{figure}[htbp]
\begin{center}
\setlength{\unitlength}{1cm}
\includegraphics[height=4cm]{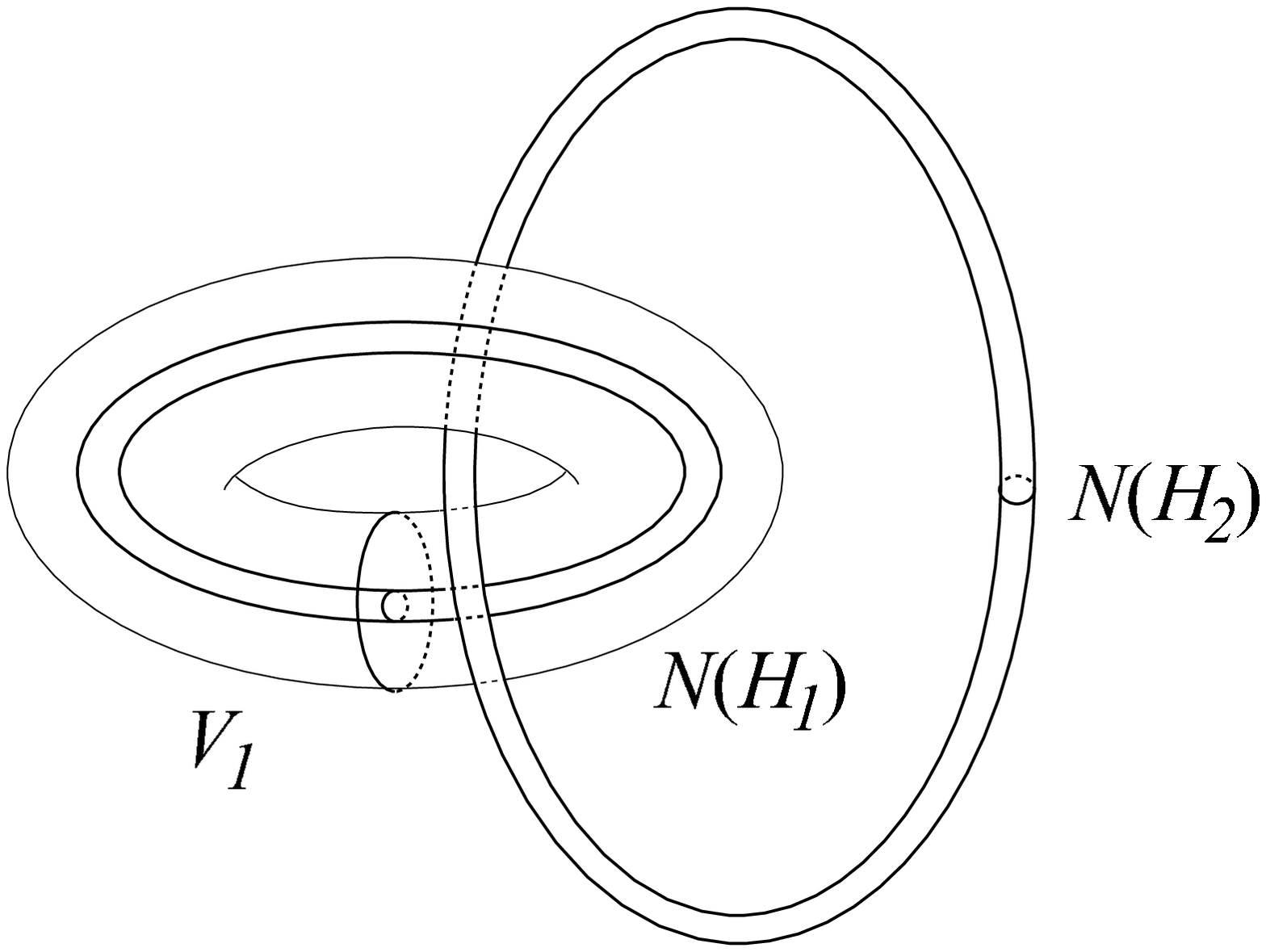}
\caption{ \label{Figure3}}
\end{center}
\end{figure}

\noindent 
{Proof}. \  
Let $h_i:D^2\times S^1\longrightarrow N(H_i)$ be the framing of $H_i$ for each $i=1,2$. Putting $j_i=h_i\vert _{S^1\times S^1}\ (i=1,2)$ and $X:=\overline{S^3-N(H_1)\cup N(H_2)}$, we have a cobordism with parametrized boundary ${\cal W}:=(X; j_1,-j_2)$. 
Then 
$$Z_{\cal W}(v_i)=\sum\limits_{i=0}^mJ(H;i,j)v_j^{\ast }.$$
 Let 
$V_1$ be the tubular neighborhood of $N(H_1)$ depicted  as in Figure \ref{Figure3}. 

Setting $V_2:=\overline{S^3-V_1}$, and defining orientation preserving diffeomorphisms $j_1':S^1\times S^1\longrightarrow \partial V_1$ and $j_2':S^1\times S^1\longrightarrow \partial V_2$ in parallel to $j_1$ and $j_2$ respectively, we see that 
\par 
$\bullet$\ $\partial V_1=-\partial V_2$,
\par 
$\bullet$\ $X=\overline{V_1-N(H_1)}\cup \overline{V_2-N(H_2)}$,
\par 
$\bullet$\ ${\cal W}_1:=(\overline{V_1-N(H_1)}; j_1, j_1')\cong Id_{S^1\times S^1}\quad \text{as cobordisms}$,
\par 
$\bullet$\ ${\cal W}_2:=(\overline{V_2-N(H_2)}; -j_2', -j_2)\cong Id_{-S^1\times S^1}\quad \text{as cobordisms}$.

\begin{figure}[htbp]
\begin{center}
\setlength{\unitlength}{1cm}
\includegraphics[width=5cm]{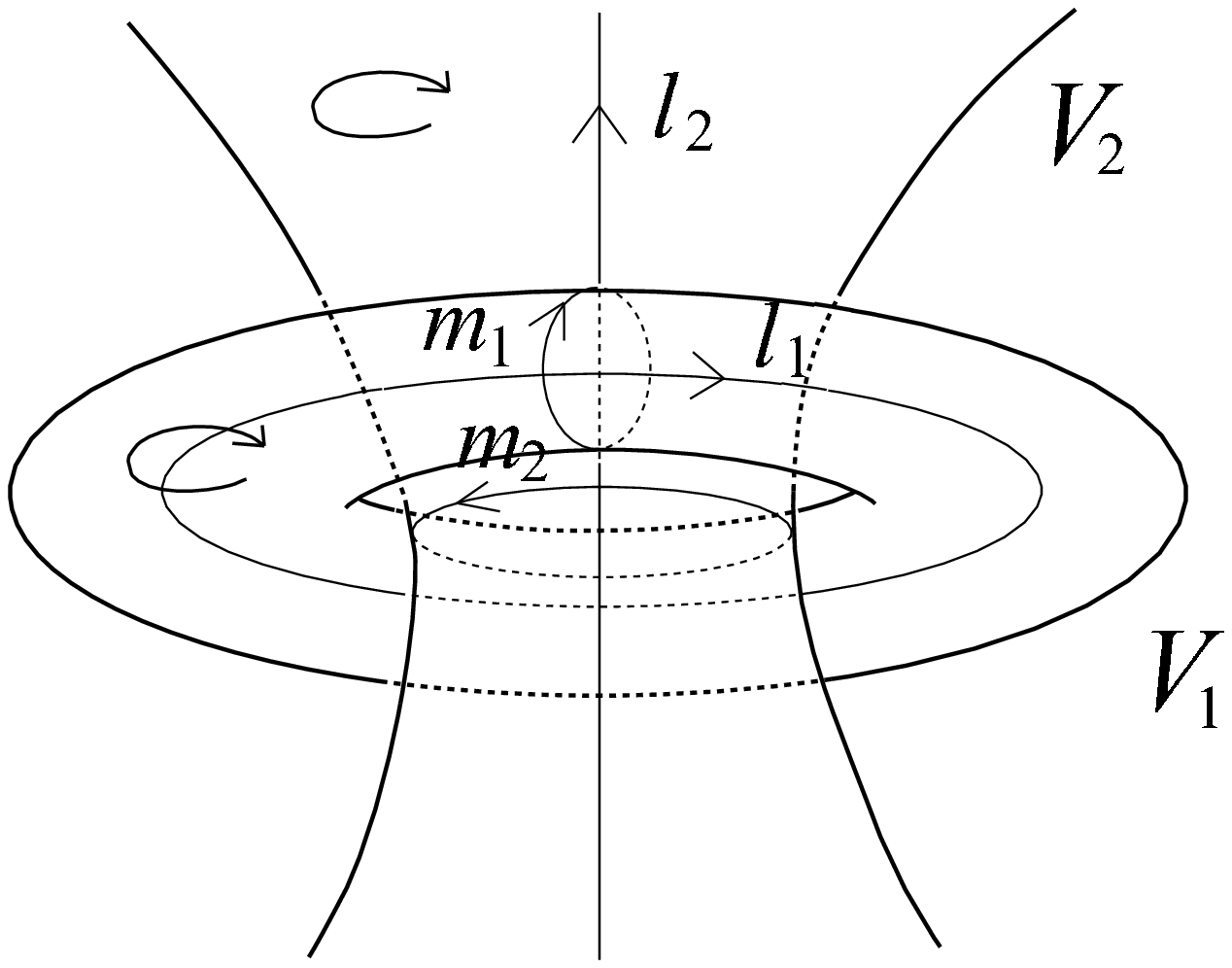}
\caption{\label{Figure4}}
\end{center}
\end{figure}

Since  
$$(-j_2')^{-1}\circ j_1':\ 
\begin{cases} 
m \longmapsto -l, \cr  
l \longmapsto -m. \end{cases}$$
for the meridian $m$ and the longitude $l$ (see Figure \ref{Figure4}), it follows that 
$$Z(-j_2'{^{-1}}\circ j_1')=Z(U\circ S).$$ 

Therefore, we obtain 
$$Z_{\cal W}=Z_{{\cal W}_2}\circ Z(U\circ S)\circ Z_{{\cal W}_1},$$
and whence 
$$\sum\limits_{j=0}^mJ(H;i,j)v_j^{\ast }=Z_{\cal W}(v_i)=(\theta \circ 
Z(U\circ S))(v_i),$$
where $\theta :Z(-S^1 \times S^1) \longrightarrow Z(S^1 \times S^1)^{\ast }$ 
denotes the natural isomorphism. 
Let $(U_{ij})_{i,j=0,1,\cdots ,m}$ be the presentation matrix of 
$\theta \circ Z(U):Z(S^1 \times S^1) \longrightarrow Z(S^1\times S^1)^{\ast }$ 
with respect to $\{ v_i \}_{i=0}^m$ and $\{ v_i^{\ast }\}_{i=0}^m$. 
Then, we have 
$$(J(H;i,j))_{i,j=0,1,\cdots ,m}=(U_{ij})(S_{ij}).$$
Hence, we have 
\begin{eqnarray*}
&(J(H;i,j))_{i,j=0,1,\cdots ,m}=(S_{ij})_{i,j=0,1,\cdots ,m} \notag \\ 
\Longleftrightarrow \ & (U_{ij})=I,\ \text{where $I$ is the identity matrix.} \notag \\ 
\Longleftrightarrow \ & (\theta \circ Z(U))(v_i)=v_{i}^{\ast } \ 
(i=0,1,\cdots ,m). \qquad \qquad \qquad \notag 
\end{eqnarray*}
This completes the proof. \hfill Q.E.D. 

\begin{figure}[hbtp]
\begin{center}
\setlength{\unitlength}{1cm}
\includegraphics[height=6cm]{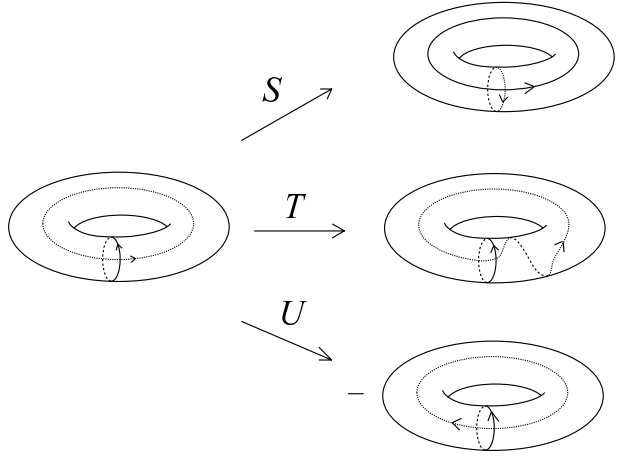}
\caption{the action of $SL_2(\mathbb{Z})$ \label{Figure6}}
\end{center}
\end{figure}

Let us recall that the mapping class group $\varGamma_{S^1\times S^1}$ of the torus $S^1\times S^1$ is isomorphic to the group 
$SL_2(\mathbb{Z})$ of integral $2\times 2$-matrices with determinant $1$. It is well-known that this group is generated by $S=\begin{pmatrix} 0 & 1 \\ -1 & 0
\end{pmatrix}$ and $T=\begin{pmatrix} 1 & 0 \\ 1 & 1\end{pmatrix}$ with relations
$S^4=I,\ (ST)^3=S^2$. The matrices $S$ and $T$ correspond to the orientation preserving diffeomorphisms $S^1\times S^1$ to $S^1\times S^1$ which are defined by $S(z,w)=(\bar{w},z)$ and $T(z,w)=(zw,w),\ (z,w)\in S^1\times S^1$, respectively, where we regard $S^1$ as the set of complex numbers of absolute value $1$ . 
To define the Verlinde basis, we need one more orientation preserving diffeomorphism $U: S^1\times S^1\longrightarrow - S^1\times S^1$ defined by $U(z,w)=(z,\bar{w})$ for $(z,w)\in S^1\times S^1$ (see Figure \ref{Figure6}). 

\par \bigskip 
\begin{defn}\label{vb}\ \ 
Let $Z$ be a $(2+1)$-dimensional TQFT. A basis $\{ v_i\} _{i=0}^m$ of $Z(S^1\times S^1)$ is said to be a {\it Verlinde basis} if it has the following properties. 
\begin{enumerate} 
\item[(i)] $v_0$ is the identity element of the fusion algebra associated 
with $Z$. 
\item[(ii)] 
\begin{enumerate} 
\item[(a)] $Z(S)$ is presented by a unitary and symmetric matrix with respect to the basis $\{ v_i\} _{i=0}^m$.
\item[(b)] $Z(S)^2v_0=v_0$, and $Z(S)^2v_i\in \{ v_j\} _{j=0}^m$ for all $i$. 
\item[(c)] We define $S_{ji}\in \mathbb{C}$ by $Z(S)v_i=\sum_{i=0}^mS_{ji}v_j$. Then, 
\begin{enumerate} 
\item[1.] $S_{i0}\not= 0$ for all $i$. 
\item[2.] $N_{ij}^k:=\sum_{l=0}^m\frac{S_{il}S_{jl}\overline{S_{lk}}}{S_{0l}}\ (i,j,k=0,1,\cdots ,m)$ coincide with the structure constants of the fusion algebra with respect to $\{ v_i\} _{i=0}^m$. 
\end{enumerate}
\end{enumerate}
\item[(iii)] $Z(T)$ is presented by a diagonal matrix with respect to the basis $\{ v_i\} _{i=0}^m$. 
\item[(iv)] $Z(U)v_i=v_{i}^{\ast }$ for all $i$ under the identification $Z(-S^1\times S^1)\cong Z(S^1\times S^1)^{\ast }$. 
\end{enumerate}
We call a {\it pre-Verlinde basis} a basis $\{ v_i\} _{i=0}^m$ of $Z(S^1\times S^1)$ satisfying the four conditions (i), (ii.a), (ii.b) and $Z(U)v_0=v_0^{\ast }$. 
\end{defn}

\par \bigskip 
\begin{rems}
1. If $\{ v_i\} _{i=0}^m$ is a pre-Verlinde basis, then $S_{i0}$ is a real number for all $i$.  
\par 
\noindent 
2. A Verlinde basis is unique up to order of elements, since $w_i=S_{0i}\sum_{j=0}^m\overline{S_{ji}}v_j$\ $(i=0,1,\cdots ,m)$ are all orthogonal primitive idempotents in the fusion algebra satisfying $1=w_0+w_1+\cdots +w_m$. This fact follows from that the $S$-matrix diagonalizes the fusion rules in conformal field theory \cite{Verlinde}. 
\par 
\noindent 
3. The map $\overline{\mathstrut \ \cdot \ }:\{ 0,1,\cdots ,m\} \longrightarrow \{ 0,1,\cdots ,m\} $ defined by $Z(S)^2v_i=v_{\bar{i}}$ is an involution satisfying $\bar{0}=0$. 
\par 
\noindent 
4. The last condition (iv) was introduced in \cite{Wakui} and modified in \cite{SuzukiWakui}. 
\end{rems}

\par \bigskip 
Let us describe some basic results on a $(2+1)$-dimensional TQFT with a 
pre-Verlinde basis.

\par \bigskip 
\begin{lem}
\label{Lemma3}
Let $Z$ be a (2+1)-dimensional TQFT and $\{ v_i\} _{i=0}^m$ a pre-Verlinde basis. 
Then, for the cobordism $W=(-D^2\times S^1; S^1\times S^1, \emptyset )$, we have  
$Z_W=v_0^{\ast }: Z(S^1\times S^1)\longrightarrow \mathbb{C}$. In particular, $Z(S^2\times S^1)=1$. 
\end{lem}

\par \medskip \noindent 
{Proof}.\ \ 
By considering the orientation preserving diffeomorphism $\tilde{U}: D^2\times S^1\longrightarrow -D^2\times S^1$ defined by $\tilde{U}(z, w)=(z,\bar{w})$, we have $(Z(\tilde{U}\vert _{\partial D^2\times S^1}))(Z(D^2\times S^1))=Z(-D^2\times S^1)\in Z(-S^1\times S^1)$. Since $\tilde{U}\vert _{\partial D^2\times S^1}=U$ and $Z(D^2\times S^1)=v_0$, it follows that $v_0^{\ast }=Z(-D^2\times S^1)=Z_W$ as elements in $Z(S^1\times S^1)^{\ast }$. 
\par 
Next, we prove that $Z(S^2\times S^1)=1$. The $3$-manifold $S^2\times S^1$ is regarded as $(D_{+}^2\times S^1)\cup (-D_{-}^2\times S^1)$, where $D_{+}^2=\{ (x,y,z)\in S^2 \ \vert \ z\geq 0\} ,\ D_{-}^2=\{ (x,y,z)\in S^2 \ \vert \ z\leq 0\} $. Thus, for the cobordisms $W_1:=(D_{+}^2\times S^1; \emptyset, \partial D_{+}^2\times S^1)$ and $W_2:=(-D_{-}^2\times S^1; \partial D_{-}^2\times S^1,\emptyset )$, we have  
$$Z(S^2\times S^1)=Z_{W_2}\circ Z_{W_1}.$$ 
Since $W_1\cong (D^2\times S^1;\emptyset ,S^1\times S^1)=W_0$ and $W_2\cong (-D^2\times S^1;S^1\times S^1,\emptyset )=W$, it follows that 
$$Z(S^2\times S^1)=Z_{W}\circ Z_{W_0}=v_0^{\ast }(v_0)=1.$$ 
This completes the proof. \hfill Q.E.D. 

\par \bigskip 
\begin{lem} 
\label{Lemma4}
Let $Z$ be a (2+1)-dimensional TQFT and $\{ v_i\} _{i=0}^m$ a pre-Verlinde basis. Then, 
$J(\text{\LARGE $\bigcirc $};i)=S_{i0}$ for each $i\in \{ 0,1,\cdots ,m\} $.
\end{lem}

\par
\begin{figure}[htbp]
\vspace{0.5cm}
\begin{center}
\setlength{\unitlength}{1cm}
\scalebox{0.9}[0.9]{\includegraphics[height=6cm]{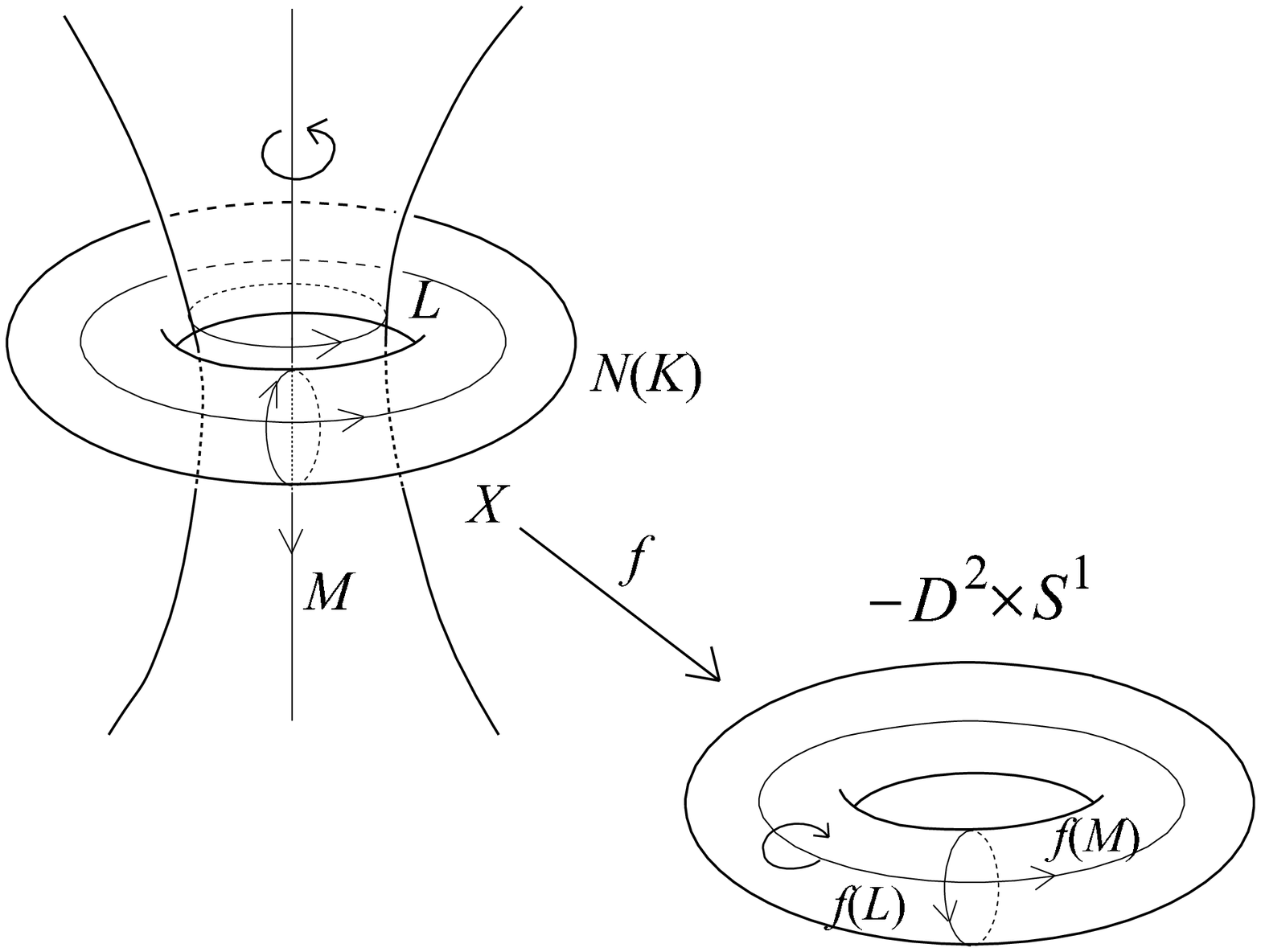}}
\caption{\label{Figure7}}
\end{center}
\end{figure}

 \medskip 
\noindent 
{Proof}. \ 
Let $K$ be the trivial knot with $0$-framing given by unit circle in $\mathbb{R}^3$, and $X=\overline{S^3-N(K)}$ the knot exterior of $K$. 
The cobordism $W_K=(X; \partial N(K), \emptyset )$ is isomorphic to the cobordism $W:=(-D^2\times S^1; S^1\times S^1, \emptyset )$ via the orientation preserving diffeomorphism $f:X\longrightarrow -D^2\times S^1$ such that  
$f(M)=l$ and $f(L)=-m$, where $M$ and $L$ are simple closed curves on $\partial N(K)$ depicted in Figure \ref{Figure7}.

\par 
Then, we have 
$$Z_{W_K}=Z_{W}\circ Z(f\vert _{\partial N(K)}).$$
Since 
$$Z(f\vert _{\partial N(K)})\circ Z(j)=Z(S)$$
for $j:=h\vert _{S^1\times S^1}$, where $h:D^2\times S^1\longrightarrow N(K)$ is the framing of $K$,   
we have 
$$J(K,i)=\langle v_0^{\ast },\ Z(S)(v_i)\rangle =S_{i0}.$$
This completes the proof.  \hfill Q.E.D. 

\par \bigskip 
The framed link invariants $J(L;i_1,\cdots ,i_r)$, $(i_1, \dots, i_r=0,1,\dots, 
m)$ have the following nice properties. 

\par 
\begin{figure}[hbtp]
\begin{center}
\setlength{\unitlength}{1cm}
\includegraphics[height=2cm]{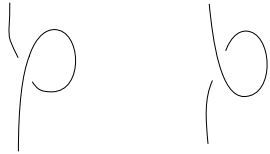} 
\caption{a positive curl and a negative curl \label{Figure8}}
\end{center}
\end{figure}

\par \bigskip 
\begin{lem} 
\label{Lemma5}
Let $Z$ be a $(2+1)$-dimensional TQFT and $\{ v_i\} _{i=0}^m$ a pre-Verlinde basis. 
For a framed link $L=L_1\cup L_2\cup \cdots \cup L_r$, the invariant $J(L; i_1, \cdots , i_r)$ has the following properties.  
\par 
(1) Let $L'=L_1\cup \cdots \cup (-L_k)\cup \cdots \cup L_r$ be the framed link obtained from $L$ by changing the orientation for the $k$-th component $L_k$. Then, 
$$J(L;i_1,\cdots ,i_k,\cdots ,i_r)=J(L';i_1, \cdots ,\bar{i}_k,\cdots ,i_r).$$
\indent 
(2) Suppose that $\{ v_i\} _{i=0}^m$ satisfies the condition (iii) in the definition of Verlinde basis. We define $t_i\in \mathbb{C}\ (i=0,1,\cdots ,r)$ by $Z(T)v_i=t_iv_i$. 
Let $L'$ be the framed link obtained from $L$ such that it is same as $L$ except for a small segment of the $k$-th component $L_k$ and the small segment is replaced by a positive or negative curl as shown in Figure \ref{Figure8}. Then, the following equations hold. 
\par 
If the small segment is replaced by a positive curl, then  
$$J(L'; i_1,\cdots ,i_r)=t_{i_k}^{-1}J(L;i_1,\cdots ,i_r).$$
If the small segment is replaced by a negative curl, then 
$$J(L'; i_1,\cdots ,i_r)=t_{i_k}J(L;i_1,\cdots ,i_r).$$
\end{lem}

\noindent 
{Proof}.\  
Without loss of generality, we may suppose that $k=1$. 

\begin{figure}[hbtp]
\begin{center}
\setlength{\unitlength}{1cm}
\includegraphics[height=5cm]{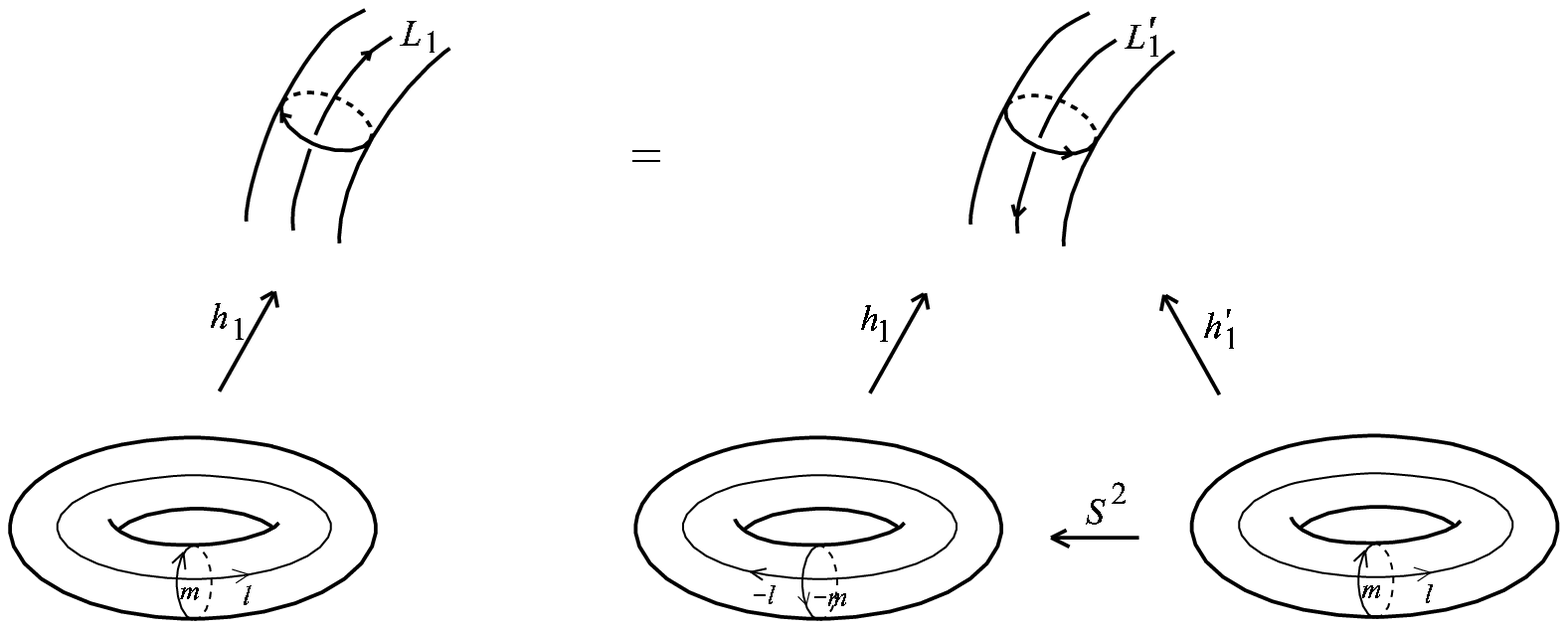}
\caption{\label{Figure9}}
\end{center}
\end{figure}

\noindent
(1) Let $h_i:D^2\times S^1\longrightarrow N(L_i)$ be the framing of $L_i$ for each $i=1,\cdots ,r$, and $h_1':D^2\times S^1\longrightarrow N(-L_1)$ the framing of $-L_1$. We set $j_i:=h_i\vert _{S^1\times S^1}$ for each $i=1,\cdots ,r$, and $j_1':=h_1'\vert _{S^1\times S^1}$. 
Then, 
$$j_1'=j_1\circ S^2,$$
since $h_1'=h_1\circ f$ for the orientation preserving diffeomorphism $f:D^2\times S^1\longrightarrow D^2\times S^1$ defined by $f(z,w)=(\bar{z},\bar{w})$ (see Figure \ref{Figure9}).

Let $X$ be the link exterior of $L$. 
Since the two cobordisms 
$${\cal W}_L=(X; \coprod_{i=1}^rj_i, \emptyset ),\quad 
{\cal W}_{L'}=(X;j_1'\coprod (\coprod _{i=2}^rj_i), \emptyset )$$ 
are isomorphic via the identity map $id_X$, we see that 
$$Z_{{\cal W}_{L'}}\circ (Z(S^2)\circ id \otimes \cdots \otimes id)=Z_{{\cal W}_L}.$$
Since $Z(S^2)v_i=v_{\bar{i}}$, we have 
$$J(L;i_1,i_2,\cdots ,i_r)=J(L';\bar{i}_1, i_2,\cdots ,i_r).$$
\par 
(2) We suppose that $L'$ arises from $L$ by replacing a small segment of the first component $L_1$ by a positive curl. 
Let $h_1':D^2\times S^1\longrightarrow N(L_1')$ be the framing of $L_1'$. We set $j_1':=h_1'\vert _{S^1\times S^1}$. 

\begin{figure}[hbtp]
\begin{center}
\setlength{\unitlength}{1cm}
\includegraphics[height=6cm]{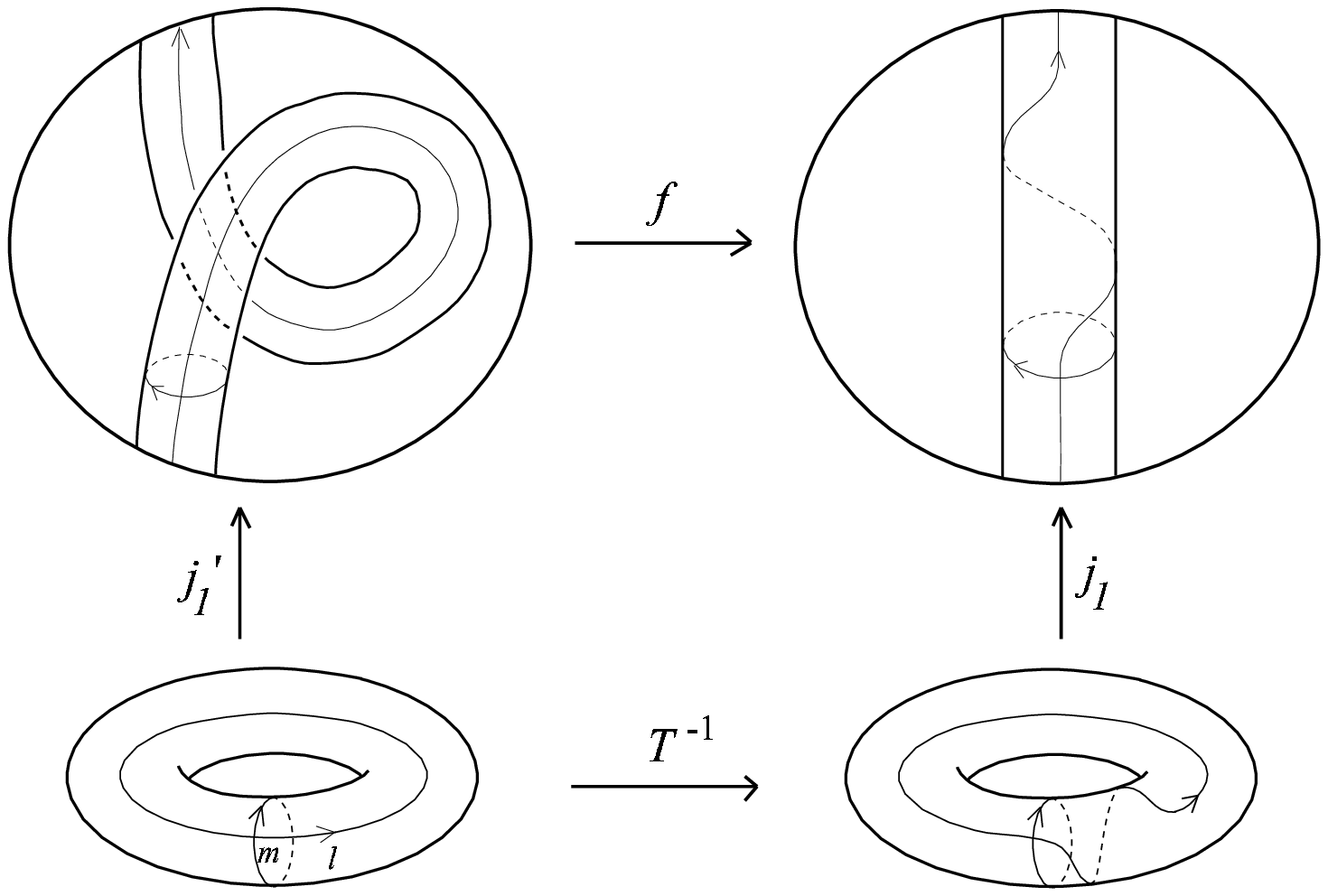} 
\caption{\label{Figure10}}
\end{center}
\end{figure}

Let $X'$ be the link exterior of $L'$.
Then, there exists an orientation preserving diffeomorphism
$f:X'\longrightarrow X$
such that
\begin{enumerate}
\item[(i)] $j_1^{-1}\circ (-f\vert _{\partial N(L_1')})\circ j_1'=T^{-1}$ up
to isotopy,
\item [(ii)] $-f\vert _{\partial N(L_i)}=id_{\partial N(L_i)}$ for all
$i=2,\cdots ,r$.
\end{enumerate}

This map $f$ gives rise to the isomorphism between
the cobordisms ${\cal W}_{L'}=(X';\emptyset ,j_1'\coprod (\coprod
_{i=2}^rj_i))$ and ${\cal W}_L=(X;\emptyset ,\coprod _{i=1}^rj_i)$. Hence,
we have
$$Z_{{\cal W}_{L'}}\circ (Z({j_1'}^{-1})\otimes (\bigotimes
_{i=2}^rZ(j_i^{-1})))=Z_{{\cal W}_L} \circ (\bigotimes
_{i=1}^rZ(j_i^{-1})\circ Z(-f\vert _{\partial N(L_i)})).$$
It follows that
$$Z_{{\cal W}_{L'}}=Z_{{\cal W}_L} \circ (Z(T^{-1})\otimes id \cdots \otimes
id ).$$
This implies that
$$J(L'; i_1,\cdots ,i_r)=t_{i_1}^{-1}J(L;i_1,\cdots ,i_r).$$

Thus, the proof of the first equation of (2) is completed.
By a similar argument, the second equation of (2) can be proved.
This completes the proof. \hfill Q.E.D. 

\begin{figure}[htpb]
\begin{center}
\setlength{\unitlength}{1cm}
\includegraphics[height=5cm]{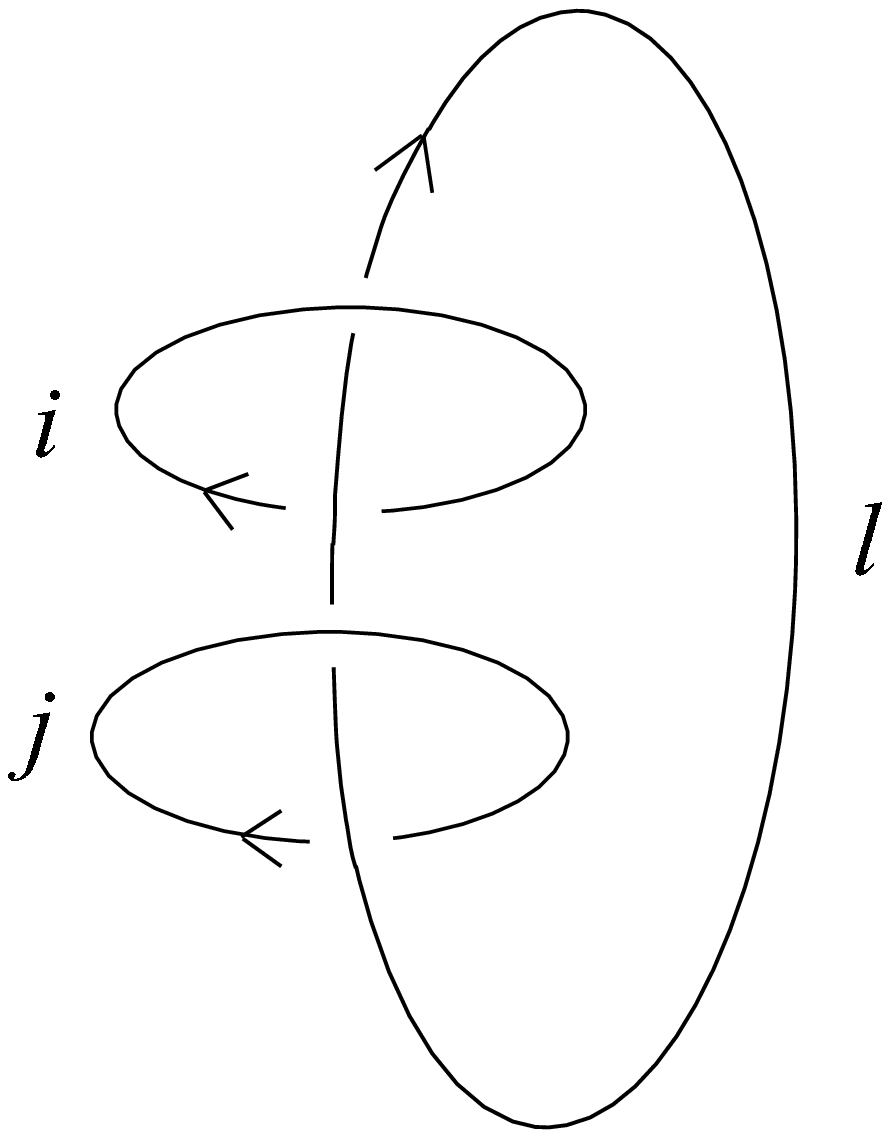}
\caption{\label{Figure11}}
\end{center}
\end{figure} 

\par \bigskip 
Let us introduce a criterion for the Verlinde identity (ii.c.2). 
\par 
We suppose that a $(2+1)$-dimensional TQFT $Z$ arises from a semisimple ribbon $\mathbb{C}$-linear Ab-category in the sense of Turaev \cite{Turaev}. Let $\{ v_i\} _{i=0}^m$ be a pre-Verlinde basis satisfying the equivalent conditions (1) and (2) in Proposition \ref{Proposition2}. We suppose $S_{i0}\not= 0$ for all $i=0,1,\cdots ,m$.   
\par
\begin{figure}[hbtp]
\begin{center}
\setlength{\unitlength}{1cm}
\scalebox{0.9}[0.9]{\includegraphics[width=3.5cm]{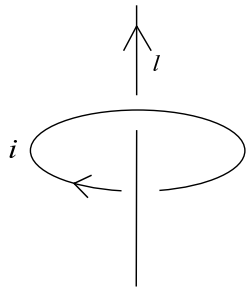}}
\caption{the colored $1-1$ tangle $T_{i,l}$ \label{Figure12}}
\end{center}
\end{figure}

 
For the $3$-component framed link $L$ presented by the diagram as in Figure \ref{Figure11}, 
the framed link invariant $J(L;i,j,k)$ coincides with the quantum trace of $J(T_{i,l})\circ J(T_{j,l})$, where  $T_{i,l}$ is the colored framed tangle presented by the diagram as in Figure \ref{Figure12}. 
\par 
Since $J(T_{i,l})$ is a map from the simple object $l$ to $l$, the invariant $J(T_{i,l})$ is a scalar multiple by $id_l$. We define $a_{i}$ to be this 
scalar. Since the quantum trace of $T_{i,l}$ is the Hopf link with colors $i,l$, we see that 
$$a_i\cdot J(\text{\LARGE $\bigcirc $};l)=S_{il}.$$
This implies that 

\begin{equation}
J(L;i,j,l)=\dfrac{S_{il}S_{jl}}{S_{l0}}. \label{eq1} 
\end{equation}

\par
\begin{figure}[htbp]
\begin{center}
\setlength{\unitlength}{1cm}
\scalebox{0.9}[0.9]{\includegraphics[height=3cm]{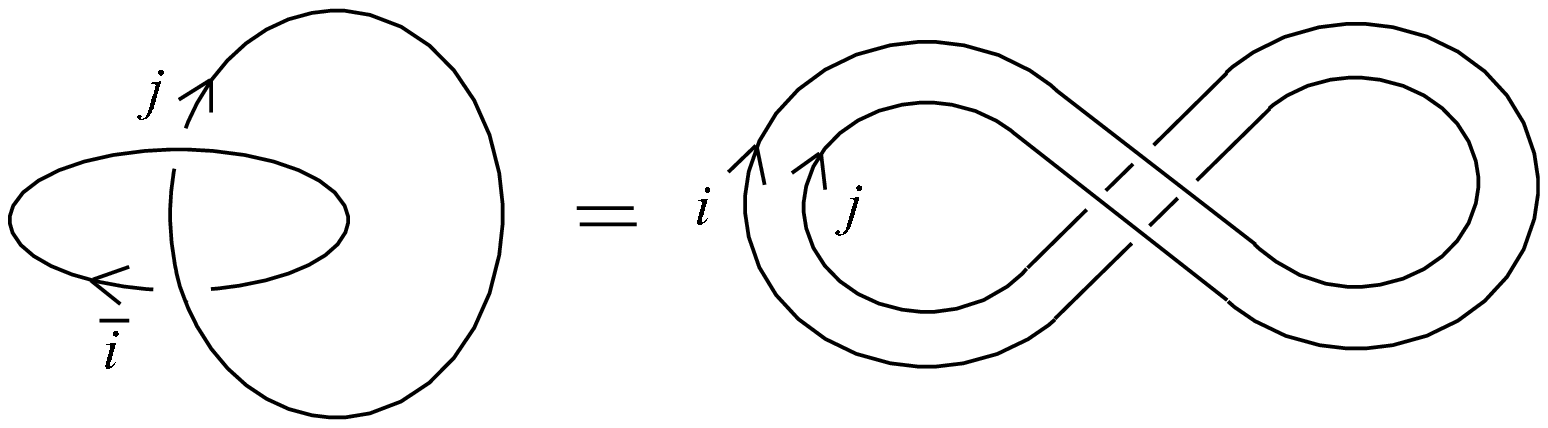}}
\caption{\label{Figure13}}
\end{center}
\end{figure}

 
On the other hand, by fusing $L_i$ and $L_j$ we have 

\begin{equation}
J(L;i,j,l)=\sum\limits_{i=0}^mN_{ij}^kS_{kl},\label{eq2} 
\end{equation}
where $N_{ij}^k\ (i,j,k=0,1,\cdots ,m)$ are the structure constants of the fusion algebra with respect to $\{ v_i\}_{i=0}^m$. 
From (\ref{eq1}) and (\ref{eq2}), it follows that 
$$\sum\limits_{i=0}^mN_{ij}^kS_{kl}=\dfrac{S_{il}S_{jl}}{S_{l0}}.$$
Since the above equation induces the Verlinde identity 
$$N_{ij}^k=\sum_{l=0}^m\frac{S_{il}S_{jl}\overline{S_{lk}}}{S_{l0}},$$
the condition (2) in Proposition \ref{Proposition2} implies the condition (ii.c.2) in the definition of Verlinde basis (see \cite{Takata}, \cite{Witten} for  similar arguments). 
If $\{ v_i\} _{i=0}^m$ satisfies the condition (iii) in the definition of Verlinde basis, then we see that the converse is true by using the technique in 
\cite{Kohno} (see Figure \ref{Figure13}). 

Thus, we conclude the following. 

\par \bigskip 
\begin{thm} 
\label{Theorem6}
Suppose that a $(2+1)$-dimensional TQFT $Z$ arises from a semisimple ribbon $\mathbb{C}$-linear Ab-category. 
Let $\{ v_i\} _{i=0}^m$ be a pre-Verlinde basis satisfying $S_{i0}\not= 0$ for all $i=0,1,\cdots ,m$. 
If $Z(U)v_i=v_{i}^{\ast }$ for all $i=0,1,\cdots ,m$, then 
$N_{ij}^k:=\sum_{l=0}^m\frac{S_{il}S_{jl}\overline{S_{lk}}}{S_{0l}}\ (i,j,k=0,1,\cdots ,m)$ coincide with the structure constants of the fusion algebra with respect to $\{ v_i\} _{i=0}^m$. 
\par 
Furthermore, if $\{ v_i\} _{i=0}^m$ satisfies the condition (iii) in the definition of Verlinde basis, then the converse is true. Therefore, the following are equivalent. 
\begin{enumerate}
\item[(i)] $J(H;i,j)=S_{ij}$ for all $i,j=0,1,\cdots ,m$. 
\item[(ii)] $Z(U)v_i=v_{i}^{\ast }$ for all $i=0,1,\cdots ,m$.
\item[(iii)] $N_{ij}^k:=\sum_{l=0}^m\frac{S_{il}S_{jl}\overline{S_{lk}}}{S_{0l}}\ (i,j,k=0,1,\cdots ,m)$ coincide with the structure constants of the fusion algebra with respect to $\{ v_i\} _{i=0}^m$. 
\end{enumerate} 
Here, $H$ is the Hopf link presented by the diagram as in Figure \ref{Figure2}. 
\end{thm}

\section{Turaev-Viro-Ocneanu $(2+1)$-dimensional TQFT (Review)}

\subsection{Sectors and finite system $\Delta$}

For a detailed exposition about sectors, see \cite{Iz}. 
Let $N \subset M$ be an inclusion of infinite factors. In this case, 
we also have a similar concept to Jones index named Kosaki 
index, and if it takes the minimal value, we write it $[M:N]_0$ as in 
the case of Jones index. In the sequel, we assmue subfactors have  
finite minimal indices. For the inclusion of $\rho(M) 
\subset M$, where $\rho \in {\rm End}(M)$, we call $[M:\rho(M)]_0^{1/2}$ 
the statistical dimension of $\rho$ and denote it by $d(\rho)$. 
We denote the set of $*$-homomorphisms from $N$ to $M$ with finite statistical 
dimensions by ${\rm Mor}(N, M)_0$.  
We say that $\rho_1$, $\rho_2 \in {\rm Mor}(N, M)_0$ 
are equivalent if there exisits a unitary $u$ in $M$ such that $u\rho_1(x) =
\rho_2(x) u$ for all $x \in N$.  This gives an equivalence relation in 
${\rm Mor}(N, M)_0$. The set of the equivalencce classes is denoted by 
${\rm Sect}(N,M)$ and its element is called an $M$-$N$ {\it sector}, which 
is denoted by $[\rho]$ for $\rho \in {\rm Mor}(N, M)_0$. In $M$-$M$ sectors 
${\rm Sect}(M,M)=:{\rm Sect}(M)$, we have the product 
$[\rho_1] \cdot [\rho_2] = [\rho_1 \circ \rho_2]$ and the summation $[\rho_1] 
\oplus [\rho_2]$. (See \cite{Iz} for the definition of the summation of 
sectors.) These operations  define a semiring structure in 
${\rm Sect}(M)$. For $[\rho_1], [\rho_2] \in {\rm Sect}(M)$, we define 
\[
 (\rho_1, \rho_2)=\{ V \in M| V\rho_1(x) = \rho_2(x) V, \forall x \in M  \}.
\]
It is called the {\it intertwiner space} between $\rho_1$ and $\rho_2$. 
If $(\rho, \rho) = \mathbb{C}1_M$, we say that $\rho$ is {\it irreducible}. 
The intertwiner space $(\rho_1, \rho_2)$ has an inner product $\langle V, W 
\rangle=W^* \cdot V$, $V, W \in (\rho_1, \rho_2)$, if $\rho_1$ is 
irreducible.  

Moreover, ${\rm Sect}(M)$ has a conjugation $\overline{[\rho]}=
[\bar{\rho}]$. Namely, for $[\rho] \in {\rm Sect}(M)$, there exist an 
endomorphism $\bar{\rho} \in {\rm End}(M)$ and a pair of intertwiners 
$R_\rho \in (id, \bar{\rho} \rho)$ and $\bar{R}_\rho \in 
(id, \rho \bar{\rho})$ such that $\bar{R}_\rho^* \rho(R_\rho)=
R_\rho^* \bar{\rho}(\bar{R}_\rho)=1/d(\rho)$. With these operations, 
${\rm Sect}(M)$ becomes a $*$-semiring over $\mathbb{C}$. 

An important thing is that ${\rm Sect}(M)$ is closed under the operations such 
as product, direct sum, irreducible decommposition and conjugation. 

Let us introduce the notion of a finite system $\Delta$ of ${\rm End}(M)_0$, 
which is a basic data to describe a topological quantum field theory from 
subfactors. Since the embedding $\iota : N \hookrightarrow M$ is an element of 
${\rm Mor}(N, M)_0$, we can consider the sector $[\iota] \in {\rm Sect}(N, M)$. 
We note that the conjugation $\overline{[\iota]}$ is an element of 
${\rm Sect}(M, N)$, and  the product $[\iota] \overline{[\iota]}$ becomes an 
element of ${\rm Sect}(M)$. 
In a similar way, $\overline{[\iota]} [\iota]$ becomes an element of 
${\rm Sect}(N)$. By decomposing  $([\iota] \overline{[\iota]})^n$, $([\iota] 
\overline{[\iota]})^n [\iota] $, $\overline{[\iota]} ([\iota] 
\overline{[\iota]})^n $ and $(\overline{[\iota]} [\iota])^n$ into irreducible 
sectors, we get $M$-$M$, $M$-$N$, $N$-$M$ and $N$-$N$ sectors 
responsibly. (A sector $[\rho]$ is said to be {\it irreducible} if 
$\rho$ is irreducible.) If the number of the irreducible sectors in the above 
decompositions is finite, then the subfactor $N \subset M$ is called of 
{\it finite depth}. Throughout this paper, we only consider subfactors with 
finite depth and finite index. For a finite depth subfactor, we 
get finitely many irrducible $M$-$M$ sectors. In other words, we have a 
representative set $\Delta=\{ \rho_\xi \}_{\xi \in \Delta_0}$ of finitely many 
irreducible $M$-$M$ sectors such that \\
(i) $[\rho_\xi]=[\rho_\eta]$ if and only if $\xi=\eta$ \\
(ii) There exists $e \in \Delta_0$ such that $\rho_e =id$ \\
(iii) For any $\xi \in \Delta_0$, there exists $\bar{\xi} \in \Delta_0$ 
such that $\overline{[\rho_\xi]}=[\rho_{\bar{\xi}}]$ \\
(iv) There exist non-negative intetgers $N_{\xi,\eta}^\zeta$ such that 
$[\rho_\xi][\rho_\eta]=\oplus_{\zeta \in \Delta_0} N_{\xi,\eta}^\zeta 
[\rho_\zeta]$. \\

\noindent
We call this $\Delta$ a finite system of ${\rm End}(M)_0$ or simply a 
{\it finite system}. 

\begin{rem}
We note that $\Delta$ can be seen as a $C^*$-tensor category in the following 
manner. The objects of the category are $\mathbb{C}$-linear span of elements of 
$\Delta$,  the morphisms of the category are intertwiners, and the 
tensor product structure is given by the compositions of endomorphisms 
in objects of the category. By an abuse of notation, we denote this 
category by $\Delta$.  See \cite{LR} for the details. 
\end{rem}
\subsection{Turaev-Viro-Ocneanu TQFT}

We need some preparations before constructing the Turaev-Viro-Ocneanu 
TQFT. In the sequel, we write $\xi$ instead of $\rho_\xi$ and so forth, 
for simplicity. 

Let $\Delta$ be the finite system of ${\rm End}(M)_0$. We consider the following 
diagram.
\[
 \begin{CD}
\xi \cdot  \alpha \cdot \eta 
@< \xi(T_1) << \xi \cdot \beta \\
@A{T_3}AA  @AA{T_2}A \\
 \gamma \cdot \eta 
@< T_4 << \delta
 \end{CD}
\]
Here, $\alpha, \beta, \gamma, \xi, \eta, \delta 
\in \Delta$ and $T_1 \in (\beta, \alpha \cdot \eta)$, 
$T_2 \in (\delta, \xi \cdot \beta)$, $T_3 \in (\gamma, 
\xi \cdot \alpha)$, $T_4 \in (\delta, \gamma \cdot \eta)$. 
Then, the composition of intertwiners $T_4^* \cdot T_3^* \cdot \xi(T_1) 
\cdot T_2$ belongs to $(\delta, \delta)$. Since $\delta$ is assumed 
to be irreducible, this composition of intertwiners is regarded as  
a complex number. We call this number a {\it quantum $6j$-symbol}.

\begin{figure}[htbp]
  \begin{center}
\scalebox{0.5}{
\includegraphics{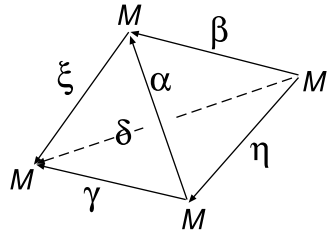}}  
 \end{center}
 \caption{A tetrahedron as a diagram}
 \label{6j}
\end{figure}

The above diagram can be seen as a tetrahedron as in Figure \ref{6j}. 

We assign $d(\beta)^{-1/2}d(\gamma)^{-1/2} T_4^* \cdot T_3^* \cdot \xi(T_1) 
\cdot T_2$ to this tetrahedron. 

Let $V$ be an oriented closed 3-dimensional manifold. Choose a triangulation 
of $V$ and write it ${\cal T}$. To each vertex in ${\cal T}$, we assign  
the factor $M$, to each edge in ${\cal T}$ an element in $\Delta$ and to each face in 
${\cal T}$ an intertwiner. 
Let $E$ be the set of the edges in ${\cal T}$, $e$ be an assignment of 
elements in $\Delta$ to the edges in ${\cal T}$, and $\varphi$ be an 
assignment of intertwiners to the faces in ${\cal T}$. 
A tetrahedron $\tau$ has the labeled edges by the assignment $e$ and 
the labeled faces by the assignment $\varphi$. Hence, to a tetrahedron 
$\tau$, we assign a complex number defined by using the $6j$-symbol as above 
or its complex conjugate depending on the orientation of $\tau$. 
We denote this complex number by $W(\tau;e,\varphi)$.  We multiply these 
$W$'s and the weights $\prod_E d(\xi)$. Then, sum up all these resulting 
values. Finally, we multiply some weights coming from the vertices 
in ${\cal T}$. 
\[
 Z_\Delta(V,{\cal T})= \lambda^{-a} \sum_e(\prod_E d(\xi))
\sum_\varphi \prod_\tau W(\tau;e,\varphi), 
\]
where $\lambda=\sum_{\xi \in \Delta_0} d(\xi)^2$ and $a$ is 
the number of the vertices in ${\cal T}$. This $Z_\Delta(V,{\cal T})$ is 
proven to be independent of any choice of triangulations because of 
the properties 
of quantum $6j$-symbols. (See \cite{EK} for a detailed account.) Namely, 
$Z_\Delta(V,{\cal T})$ turns out to be a topological invariant of $V$. So, we 
drop ${\cal T}$ off from $Z_\Delta(V,{\cal T})$ and denote this value 
by $Z_\Delta(V)$. We call $Z_\Delta(V)$ the {\it Turaev-Viro-Ocneanu invariant} 
of the 3-dimensional manifold $V$. 
When $V$ is an oriented, compact 3-dimensional manifold possibly with 
boundary, first we fix a triangulation  $\partial {\cal T}$ of the 
boundary of $V$ and extend it to the whole triangulation of $V$. Then, 
we assign an element in $\Delta$ to each edge in $\partial {\cal T}$ and 
assign an intrertwiner to each face. We fix these assignments to the end 
and denote them by $\partial e$ and $\partial \varphi$, respectively. In a similar 
way to the closed case, we assign $M$, an element in $\Delta$ and an 
intertwiner to each vertex, each edge and each face in the triangulation 
${\cal T} \setminus \partial {\cal T}$. We make 
$Z_\Delta(V,{\cal T}, \partial e, \partial \varphi)$ as above:
\[
 Z_\Delta(V,{\cal T}, \partial e, \partial \varphi)= 
 \lambda^{-a+\partial a/2} \prod_{\partial e} 
 d(\xi)^{1/2} \sum_{e \setminus \partial e}(\prod_{E \setminus \partial E} 
d(\xi)) \sum_\varphi  \prod_\tau W(\tau;e,\varphi), 
\]
where $\partial a$ is the number of the vertices on the boundary. This 
value $Z_\Delta(V,{\cal T}, \partial e, \partial \varphi)$ does not 
depend on the assignments of the vertices, the edges and the faces in 
${\cal T} \setminus \partial {\cal T}$. Such extended 
$Z_\Delta$ gives rise to a unitary TQFT because to an oriented closed surface, 
it assigns a finite dimensional Hilbert space with the inner product 
induced from the space of intertwiners. (See \cite{EK} for the detailed 
construction.) We call this TQFT the {\it $(2+1)$-dimensional 
Turaev-Viro-Ocneanu TQFT} and denote it by $Z_\Delta$, again.

\section{Verlinde basis of $Z_{\Delta}(S^1 \times S^1)$}

Let $N \subset M$ be a subfactor of an infinite factor $M$ with finite 
index and finite depth, and let $\Delta$ be a finite system of irreducible 
$M$-$M$ endomorphisms arising from the subfactor. We write, for instance, 
$\xi$ instead of $\rho_\xi$ and so forth for the elements of $\Delta$. 

Based on $\Delta$, we construct a new finite dimensional $C^*$-algebra named 
the {\it tube algebra} ${\rm Tube} \Delta$ as in \cite{Ocneanu}. (Also see 
\cite{Izumi1}, but our normalization convention is different from that 
there.)  In this section, ${\rm Tube} \Delta$ plays a crutial role to find a 
nicely behaved basis of $Z_{\Delta}(S^1 \times S^1)$, which we call a Verlinde 
basis. It makes the Turaev-Viro-Ocneanu TQFT a rich theory. 

We sometimes use the simple notation $Z(S^1 \times S^1)$ instead of 
$Z_{\Delta}(S^1 \times S^1)$ in the sequel.

\subsection{Tube algebras}

A tube algebra ${\rm Tube}\Delta$, which was first introduced by Ocneanu 
in \cite{Ocneanu}, is defined by $\bigoplus_{\xi,\eta,\zeta} (\xi \cdot 
\zeta, \zeta \cdot \eta)$ as a vector space over $\mathbb{C}$. Its element 
is a linear combination of the composition  $T_2 \cdot T_1^*$ of the orthonormal 
bases of intertwiner spaces $T_1 \in (\delta, \xi \cdot \zeta)$, $T_2 \in (\delta, 
\zeta \cdot \eta)$, where $\xi, \eta, \zeta$ and $\delta$ run over $\Delta_0$. 
An element in ${\rm Tube}\Delta$ can be depicted as in the left-hand side of 
Figure \ref{tube}. We will define a product structure and a $*$-structure on it.

\begin{figure}[htbp]
  \begin{center}
\scalebox{0.35}{
\includegraphics{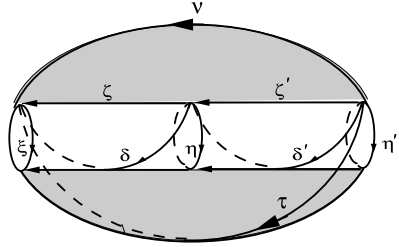}}  
 \end{center}
 \caption{A coeffcient of the product of $X \cdot Y$}
 \label{tube-c}
\end{figure}

\begin{figure}[htbp]
  \begin{center}
\scalebox{0.5}{
\includegraphics{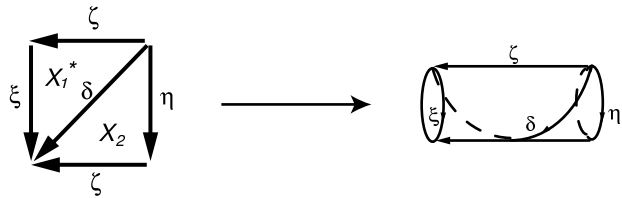}}  
 \end{center}
 \caption{An element of the tube algebra as a tube}
 \label{tube}
\end{figure}

\noindent
\underline{The product structure} \\
Let $X=X_2 X_1^* \in (\xi \cdot \zeta, \zeta \cdot \eta)$ for $X_1 \in 
(\delta, \xi \cdot \zeta)$, $X_2 \in (\delta, \zeta \cdot \eta)$ and 
$Y=Y_2 Y_1^* \in (\xi' \cdot \zeta', \zeta' \cdot \eta')$ for $Y_1 \in 
(\delta', \xi' \cdot \zeta')$, $Y_2 \in (\delta', \zeta' \cdot \eta')$. 
Then, the product of $X$ and $Y$ in ${\rm Tube} \Delta$ is defined by 
the following formula. 
\[
 X \cdot Y = \delta_{\eta, \xi'} \lambda \sum
d(\xi)^{-1/2} d(\eta')^{-1/2}  d(\eta) \lambda(X,Y;Z)Z.
\] 
Here, the summation is taken over $Z=Z_2 Z_1^*$, $Z_1 \in (\tau, \xi \cdot \nu)$, 
$Z_2 \in (\tau, \nu \cdot \eta')$, where $Z_1$ and $Z_2$ are orthonormal bases of 
the intertwiner spaces $(\tau, \xi \cdot \nu)$ and $(\tau, \nu \cdot 
\eta')$ respectively, and $\lambda(X,Y;Z)$ is the value of $Z_\Delta$ of 
the 3-manifold depicted in Figure \ref{tube-c}. 
For general elements $X$, $Y$ in ${\rm Tube} \; \Delta$, we 
define the product of them by linearlity since $X$ and $Y$ are linear 
combinations of the forms $X_2 X_1^*$ and $Y_2 Y_1^*$ as above. 

\begin{figure}[htbp]
  \begin{center}
\scalebox{0.4}{
\includegraphics{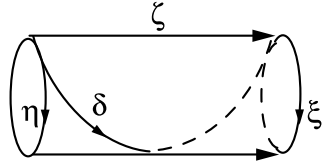}}  
 \end{center}
 \caption{An inversed tube}
 \label{tube-r}
\end{figure}

\noindent
\underline{The $*$-structure} \\
Let $X$ be as above. Then, we can consider $X$ as a tube as in the right-hand 
side of Figure \ref{tube}. 
The $*$-operation is defined by the inversing the tubes inside out. We denote 
this $*$-operation by $X^*$. See Figure \ref{tube-r}.

With the product and the $*$-structure defined as above, ${\rm Tube} \Delta$ 
becomes a finite dimensional $C^*$-algebra. Since any finite dimensional 
$C^*$-algebra is semisimple, we may assume that 
${\rm Tube} \Delta \cong \bigoplus_{i=0}^r M_{n_i}({\mathbb C})$. 

We observe that the definition of tube algebras is compatible 
with the operations of Turaev-Viro-Ocneanu TQFT such as gluing, 
cutting and so forth.

\begin{rem}
 In \cite{Izumi1}, M. Izumi has introduced the tube algebra in the 
 setting of sectors, but it is slightly different from ours in the 
normalization coefficients. We followed the definition of the tube 
algebra in \cite[Chapter 12]{EK}. 
\end{rem}

Before we start the analysis of tube algebras, we list the notations that 
we will use frequently. 

For the solid torus $D^2 \times S^1$, we denote 
the value of $Z_\Delta$ of $D^2 \times S^1$ assigned a vector $\lambda$ 
on the boundary by $Z_\Delta(D^2 \times S^1;\lambda)$. 

\begin{figure}[htbp]
  \begin{center}
\scalebox{0.4}{
\includegraphics{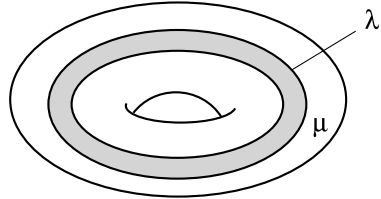}}  
 \end{center}
 \caption{The labeled 3-manifold obtained by removing another solid torus from a solid torus}
 \label{solid-a}
\end{figure}

For the 3-manifold $D^2 \times S^1 \setminus {\rm Int}D_0^2 \times S^1$, 
where $D_0$ is contained in ${\rm Int} D^2$, we denote the value of $Z_\Delta$ of 
$D^2 \times S^1 \setminus {\rm Int}D_0^2 \times S^1$ assigned vectors $\lambda$ and 
$\mu$, $\lambda$ for the boundary of $D_0^2 \times S^1$ and $\mu$ for the boundary of 
$D^2 \times S^1$ by $Z_\Delta(D^2 \times S^1 \setminus {\rm Int}D_0^2 \times S^1; 
\lambda,\mu)$. See Figure \ref{solid-a}. 

\begin{figure}[htbp]
  \begin{center}
\scalebox{0.4}{
\includegraphics{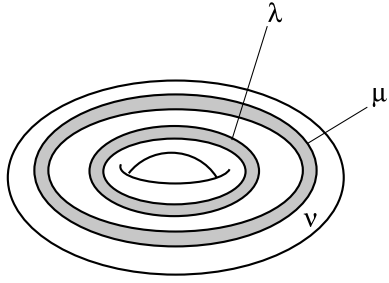}}  
 \end{center}
 \caption{The labeled 3-manifold obtained by removing two solid tori 
 from a solid torus}
 \label{solid-b}
\end{figure}

In a similar manner, for the 3-manifold $D^2 \times S^1 \setminus ( {\rm Int} D_1^2 
\times S^1 \cup {\rm Int} D_2^2 \times S^1)$, where $D_1^2$ and $D_2^2$ are two disjoint 
disks contained in ${\rm Int} D^2$, we denote the value of $Z_\Delta$ of $D^2 \times S^1 
\setminus ({\rm Int} D_1^2 \times S^1 \cup {\rm Int} D_2^2 \times S^1)$ 
assigned vectors $\lambda$, $\mu$ and $\nu$, $\lambda$ for the boundary of 
$D_1^2 \times S^1$, 
$\mu$ for the boundary of $D_2^2 \times S^1$ and $\nu$ for the the boundary of 
$D^2 \times S^1$, by $Z_\Delta(D^2 \times S^1 \setminus ({\rm Int}D_1^2 \times 
S^1 \cup {\rm Int}D_2^2 \times S^1); \lambda,\mu;\nu)$. 
See Figure \ref{solid-b}. 

\bigskip

Let $\{ \pi_0, \cdots, \pi_r \}$ be the minimal central projections of 
${\rm Tube} \Delta$. Then, we have the following lemma. 
\begin{lem}\label{innerproduct}
$\langle \pi_i, \pi_j \rangle_{Z(S^1 \times S^1)} = \delta_{i,j} n_i^2$ 
($i,j=0, \dots, r$), where $\langle \cdot, \cdot \rangle_{Z(S^1 \times S^1)}$ is 
the inner product of $Z_\Delta(S^1 \times S^1)$ defined by the Turaev-Viro-Ocneanu 
TQFT. 
\end{lem}

\noindent
Proof. First, we consider the case $i=j$. We note that $\langle \pi_i, \pi_i 
\rangle_{Z(S^1 \times S^1)}=Z_\Delta(D^2 \times S^1 \setminus {\rm Int}D_0^2 
\times S^1;\pi_i, \pi_i)$. 

Cut $D^2 \times S^1 \setminus {\rm Int}D_0^2 \times S^1$ along the meridian, 
then $T ={\rm annulus} \times [0,1]$ is created. Let $A_{jk}$ be the value of 
$Z_\Delta(T, \xi_j, \xi_k$), where $T$ is labeled by $\xi_j$ on one side of 
the sections and by $\xi_k$ on the other side. Then, $Z_\Delta(D^2 \times S^1 
\setminus {\rm Int}D_0^2 \times S^1;\pi_i, \pi_i)=\sum_{\xi_j} A_{jj}$. 
There exists an operator $A$ such that $A_{jk}=\langle A 
\xi_j, \xi_k \rangle$, where $\langle \cdot, \cdot \rangle$ is the inner 
product on intertwiners on the sections of $T$ so that $\langle \xi_j, \xi_k 
\rangle=\delta_{jk}$. 

It is easy to see that the operator $A$ is a projection. Hence, 
$\langle \pi_i, \pi_i \rangle_{Z(S^1 \times S^1)}=\sum_{j=1}^m A_{jj} =
{\rm Tr}(A) \in {\mathbb N}$. Namely, we have proved that 
$\langle \pi_i, \pi_i \rangle_{Z(S^1 \times S^1)}={\rm dim}\; \pi_i({\rm Tube} 
\Delta) = n_i^2$. 

When $i \ne j$, it is easy to see that $\langle \pi_i, \pi_j 
\rangle_{Z(S^1 \times S^1)}=0$, because $\pi_i$ and $\pi_j$ are central 
projections orthogonal to each other. This ends the proof. \hfill Q.E.D.

\begin{thm}\label{centerTube}
 Let $\Delta$ be a finite system. Then, 
 the center of ${\rm Tube}\Delta$ is naturally isomorphic to 
 $Z_\Delta(S^1 \times S^1)$ as a vector space.
\end{thm}

\noindent
Proof.  Let $V(S^1 \times S^1)$ be the linear span of elements in ${\rm Tube}\; \Delta$ such 
that the left and right labels in Figure \ref{tube} are equal. Since ${\rm Tube}\Delta 
\supset V(S^1 \times S^1)$ and $V(S^1 \times S^1)$ contains the center of 
${\rm Tube}\Delta$, it is enough to consider $V(S^1 \times S^1)$ instead of whole 
tube algebra. 

For $X \in V(S^1 \times S^1) $, we set  $\varphi_i(X)=Z_\Delta(D^2 \times S^1 \setminus 
{\rm Int}D_0^2 \times S^1; \pi_i, X)$. If $X \perp V(S^1 \times S^1)$ 
with respect to the colored inner product of ${\rm Tube} \Delta$, we set 
$\varphi_i(X)=0$. Here, for $X=X_2 X_1^* \in (\xi \cdot \zeta, \zeta \cdot \eta)$, $X_1 \in 
(\delta, \xi \cdot \zeta)$, $X_2 \in (\delta, \zeta \cdot \eta)$ and $Y=Y_2 Y_1^* \in 
(\xi' \cdot \zeta', \zeta' \cdot \eta')$, $Y_1 \in (\delta', \xi' \cdot \zeta')$, 
$Y_2 \in (\delta', \zeta' \cdot \eta')$,  the colored inner product 
$\langle X, Y \rangle_{\rm color}$ of ${\rm Tube} \Delta$ is defined by 
$\delta_{\xi,\xi'}\delta_{\eta,\eta'}\delta_{\zeta,\zeta'} \langle X_1, Y_1 \rangle
\langle X_2, Y_2 \rangle$. The last two brackets stand for the inner 
products of the intertwiner spaces $(\delta, \xi \cdot \zeta)$ and 
$(\delta, \zeta \cdot \eta)$, respectively. Then, this $\varphi_i$ is a linear functional 
defined on ${\rm Tube}\Delta$. It is easy to see that $\varphi_i$ is tracial. 

If $X$ is an element in $\pi_i({\rm Tube}\Delta)$, then we have
\[
 \varphi_i(X)=\langle \pi_i, \pi_i \rangle_{Z(S^1 \times S^1)} {\rm  tr}(X), 
\]
where ${\rm tr}$ is the normalized trace. In a similar way, for $X,Y \in 
\pi_i({\rm Tube}\Delta)$, the value 
$Z_\Delta(D^2 \times S^1 \setminus {\rm Int}D_0^2 \times S^1; Y, X)$ is 
equal to 
$\langle \pi_i, \pi_i \rangle_{Z(S^1 \times S^1)} {\rm tr}(X){\rm tr}(Y^*)$. 

When we write $X_j=\pi_jX$ for $X \in V(S^1 \times S^1)$, $X=\oplus_{j=0}^r X_j$.
For $X \in {\rm Tube}\Delta$,  we put 
\[
 E(X)= \sum_{i=0}^r \frac{\varphi_i(X)}{\langle \pi_i, \pi_i 
 \rangle_{Z(S^1 \times S^1)} } \pi_i=
 \sum_{i=0}^r {\rm tr}(X_i) \pi_i. 
\]
Then, this is a conditional expectation from ${\rm Tube}\Delta$ to the 
center of ${\rm Tube}\Delta$.  We have a description of the kernel of $E$ by 
${\rm Ker}E = {\rm the \; linear \; span \; of \; }
\{ X \in V(S^1 \times S^1)|{\rm tr}(X_j)=0, j=0, \dots, r 
\} \cup \{ X \in {\rm Tube}\Delta | \langle X, Y \rangle_{\rm color}=0 \; {\rm for \  
all}\; Y \in  V(S^1 \times S^1) \}$. 

We have 
\begin{eqnarray*}
 \langle X_j, X_j \rangle_{Z(S^1 \times S^1)}
 &=&\langle \pi_j, \pi_j 
 \rangle_{Z(S^1 \times S^1)} {\rm tr}(X_j){\rm tr}(X_j^*) \\
 &=& \langle \pi_j, \pi_j \rangle_{Z(S^1 \times S^1)} |{\rm tr}(X_j)|^2, 
\end{eqnarray*}
and we put $Q=\{  X \in V(S^1 \times S^1)|\langle X, X \rangle_{Z(S^1 \times 
S^1)}=0 \}$.  
Then, 
$Z(S^1 \times S^1)=V(S^1 \times S^1)/Q = {\rm Tube}\Delta /{\rm Ker}E
\cong {\rm Center}({\rm Tube}\Delta)$. 
Precisely, denoting the embedding map from $Z_\Delta(S^1 \times S^1)$ into 
Tube $\Delta$ by $\iota$, we have proved that Center(Tube $\Delta$) 
$=\iota(Z_\Delta(S^1 \times S^1)) \subset {\rm Tube} \; \Delta$. 
\hfill Q.E.D.

\begin{rem}
 From this theorem and Lemma \ref{innerproduct}, 
 $\{ \frac{\pi_i}{n_i} \}_{i=0}^r$ is an orthonormal basis of 
 $Z(S^1 \times S^1)$. 
\end{rem}

\subsection{Verlinde basis}

In this subsection, we show the existence of a basis of 
$Z_\Delta(S^1 \times S^1)$ nicely behaved under the action of 
$SL_2(\mathbb{Z})$, which we call the {\it Verlinde basis}. It 
deeply depends on the structure of the tube algebra. 

Before the proof of the existence of such basis, we need some preparations.

Let $p_i$ be a minimal projection in $\pi_i ({\rm Tube} \Delta)$ for 
each $i=0,\dots,r$. By Lemma \ref{innerproduct}, we have 
$\langle p_i, p_j \rangle_{Z(S^1 \times S^1)} = \delta_{ij}$. 
(We use the same notation $p_i$ as an element of $Z_\Delta(S^1 \times S^1)$.) 
Then, by the last remark in the previous subsection, we have 
\begin{equation}\label{p_i-pi_i}
 p_i=\frac{\pi_i}{n_i}\ (i=0,\dots, r)
\end{equation}
in $Z_\Delta(S^1 \times S^1)$.

\begin{figure}[htbp]
  \begin{center}
\scalebox{0.4}{
\includegraphics{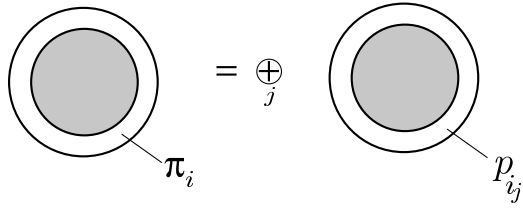}}  
 \end{center}
 \caption{The decomposition of a labeled disk}
 \label{label-pi}
\end{figure}

Let us compute the value $Z_\Delta(D^2 \times S^1; p_i)$. For this, 
we look at $Z_\Delta(D^2 \times S^1; \pi_i)$. It is a summation of 
$Z_\Delta(D^2 \times [0,1], \xi_k,\xi_k)$ over the intertwiners $\xi_k$'s 
since $D^2 \times [0,1]$ is created by cutting the solid torus along the 
meridian, This is nothing but the summation of the dimension 
of the partially labeled disks labeled by $p_{i_j}$ over $j$, $j=1,\cdots, 
n_i$, where $p_{i_j}$'s are minimal projections such that $\sum_j p_{i_j}=
\pi_i$. (See \cite[Chapter 12]{EK} for the definition of the partially labeled 
surfaces.) See Figure \ref{label-pi}. 

When we denote the dimension of the partially labeled disk with label 
$p_{i_j}$ by $c_i$, the dimension of the partially labeled disk with label 
$\pi_i$ becomes $n_i c_i$. 

\begin{figure}[htbp]
  \begin{center}
\scalebox{0.4}{
\includegraphics{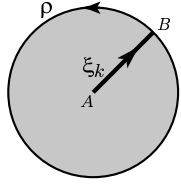}}  
 \end{center}
 \caption{A triangulated disk}
 \label{solid-d}
\end{figure}

Take a summation of $n_i c_i$ over all $i$'s, then it is equal to the 
dimension of the partially labeled disk with label $\sum_i \pi_i=1$. Hence, 
it is equal to the dimension of the triangulated disk in Figure 
\ref{solid-d}, which edge $AB$ is glued together. 
(The boudary element is a direct sum taken over arbitray $\rho$.)

Then, this vector space is non-trivial only in the case of $\rho=id$ and then, 
the dimension is one. It means that  $\sum_i n_i c_i =1$. Hence, 
only one summand can survive. We may and do assume $n_0 c_0=1$. (Hence, 
$n_0=1$, $c_0=1$.) In a similar manner, we denote $\pi_i$ by $\pi_0$ in this 
case. 

Let us summerize the above argument as a lemma. 

\begin{lem}\label{value-solid}
$Z_\Delta(D^2 \times S^1 \setminus {\rm Int} D_0^2 \times S^1;\pi_i)
=Z_\Delta(D^2 \times S^1 \setminus {\rm Int} D_0^2 \times S^1;p_i)
=\delta_{i0}$.
\end{lem}

\begin{figure}[htbp]
  \begin{center}
\scalebox{0.4}{
\includegraphics{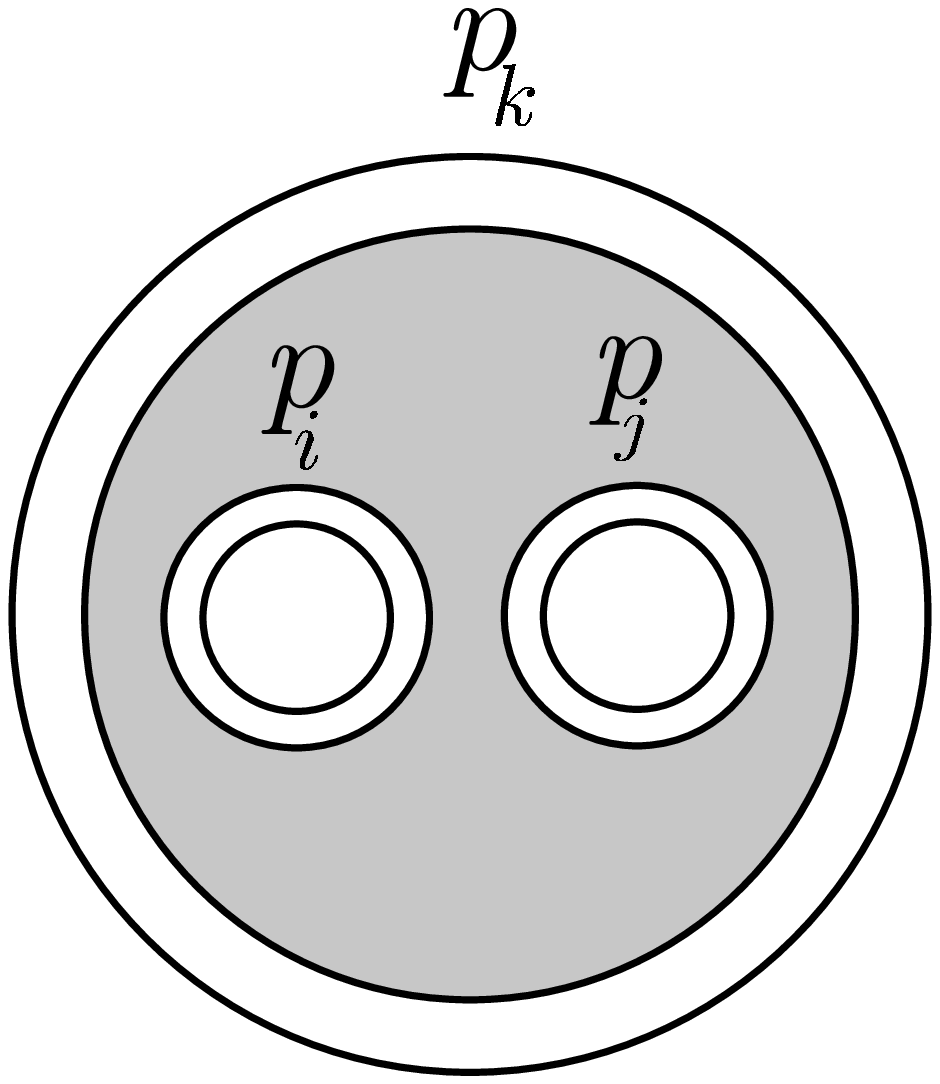}}  
 \end{center}
 \caption{$H_{\rm pants}$}
 \label{Nijk}
\end{figure}

Let $N_{ij}^k$ be $Z_\Delta (D^2 \times S^1 \setminus ( {\rm Int}D_1^2 \times 
S^1 \cup {\rm Int}D_2^2 \times S^1);p_i,p_j;p_k )$. See Figure \ref{solid-b}.  
Cut $D^2 \times S^1 \setminus ( {\rm Int}D_1^2 \times S^1 \cup 
{\rm Int}D_2^2 \times S^1)$ along the meridian, then we have $P=(D^2 \setminus 
{\rm Int} D_1^2 \cup {\rm Int} D_2^2) \times [0,1]$, and as a value, $N_{ij}^k$ is 
a summation of $Z_\Delta(P, \xi_k, \xi_k)$ over $\xi_k$'s, where $P$ is labeled 
by $\xi_k$'s on the sections. This value can be written 
$\langle A \xi_k, \xi_k \rangle$, where $\langle \cdot, \cdot \rangle$ is the inner 
product of the Hilbert space of the section. Then, it is easy to see that $A$ is a 
projection. Hence, $N_{ij}^k={\rm Tr}(A)  = {\rm dim } \; H_{\rm pants} \in {\mathbb N}$, 
where $H_{\rm pants}$ is the Hilbert space associated with the 3-holed sphere in the 
Turaev-Viro-Ocneanu TQFT. See Figure \ref{Nijk}. 

With these settings, we can prove the following theorem, which is one of 
our main theorems in this paper. (Many of the contents have already appeared 
in \cite{EK}, \cite{EK1} following several presentations of Ocneanu, but 
some normalizations are missing or incorrect there, so we include a 
complete proof here.)
\begin{thm}\label{verlinde-basis}
 Let $\Delta$ be a finite system and $\{ p_i \}_{i=0}^r$ minimal 
 projections such that each $p_i$ belongs to $\pi_i ({\rm   Tube}(\Delta))$,  
 $i=0,\dots, r$. Then, 
 $\{ p_i \}_{i=0}^r$ is a Verlinde basis of $Z_\Delta(S^1 \times S^1)$ 
 in the sense of Section 2. 
\end{thm}

\noindent
Proof. We check all the conditions for $\{ p_0, \dots, p_r \}$ to be a 
Verlinde basis step by step. 

\begin{figure}[htbp]
  \begin{center}
\scalebox{0.8}{
\includegraphics{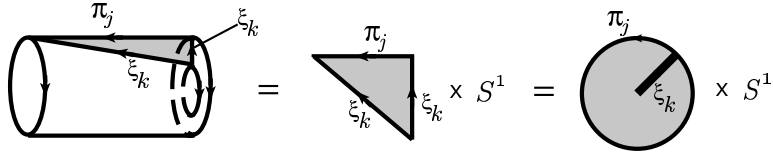}}  
 \end{center}
 \caption{$X_0$ after changes}
 \label{solid-e}
\end{figure}

Condition (ii.c.1): 
Look at the inner product $n_j^2=\langle \pi_j, \pi_j 
\rangle_{Z(S^1 \times S^1)}$. This value is defined by 
$Z_\Delta(D^2 \times S^1 \setminus {\rm Int}D_0^2 \times S^1;\pi_j,\pi_j)$. 
Cut $D^2 \times S^1 \setminus {\rm Int}D_0^2 \times S^1$ along the meridian, 
then we have $X_0={\rm annulus} \times [0,1]$ as topological object. 
Since $\pi_j$ is central, $X_0$ can be depicted as in Figure \ref{solid-e}, 
where $\xi_k$'s are intertwiners on the sections. 

Let $X$ be the solid torus that we had at last in Figure \ref{solid-e}. 
Then, we have 
\begin{equation}
 \langle \pi_j, \pi_j \rangle_{Z(S^1 \times S^1)}=\sum_{\xi_k} \alpha(\xi_k)
 Z_\Delta (X),
\end{equation}
where $\alpha(\xi_k)'s$ are positive coefficients determined by 
some products of statistical dimensions. 
It is obvious that each $Z_\Delta (X)$ is non-negative. 

Since $S(p_j)=\sum_i S_{ij} p_i$, it follows from Lemma \ref{value-solid} that 
$Z_\Delta (D^2\times S^1;S(p_j))= S_{0j}$. 

It is now easy to see  that the value of the solid torus labeled by 
$S(p_j)$ has a similar expression of the summation when it is cut along the 
meridian. Namely, $S_{0j}$ is written in the form $S_{0j}=\sum_{\xi_k} 
\beta(\xi_k) Z_\Delta (X)$, where $\beta(\xi_k)$'s are  strictly positive 
coefficients detemined by $n_j$ and some products of statistical dimensions. 
As we saw, not all of the values of $Z_\Delta (X)$ can be zero, so 
one of them must be strictly positive. Hence, $S_{0j} > 0$. 

Condition (ii.a), (iii): Since both $S$ and $T$-matrices are 
unitary, we check that $S$ is symmetric and $T$ is diagonal with repsect to 
$\{p_0, \dots, p_r \}$. 

First we prove that $S$ is symmetric. From a topological observation of 
the $Z(S)$-action on $Z(S^1 \times S^1)$, 
we have 
\begin{eqnarray*}
S_{ij}
&=&\langle S(p_i), p_j \rangle_{Z(S^1 \times S^1)} 
= \overline{\langle S(p_i)^*, p_j^* \rangle}_{Z(S^1 \times S^1)} 
= \overline{\langle S(p_i)^*, p_j \rangle}_{Z(S^1 \times S^1)} \\
&=& \overline{\langle S^*(p_i), p_j \rangle}_{Z(S^1 \times S^1)} 
= \overline{\langle p_i, S(p_j) \rangle}_{Z(S^1 \times S^1)} 
= \langle S(p_j), p_i \rangle_{Z(S^1 \times S^1)} \\
&=&S_{ji}. 
\end{eqnarray*}

Next, we prove that $T$ is diagonal. From a topological observation of 
the $Z(T)$-action on $Z(S^1 \times S^1)$, we have 
\[
 \langle T(X) \cdot T^*(Y), Z \rangle_{Z(S^1 \times S^1)} = 
\langle X \cdot Y, Z \rangle_{Z(S^1 \times S^1)}
\] 
for any $X,Y,Z \in V(S^1 \times S^1)$, where $\cdot$ stands for the multiplication 
in the tube algebra. Hence, we have $T(X) \cdot T^*(Y)=X \cdot Y$ for any 
$X, Y \in V(S^1 \times S^1)$. This implies that $T$ is diagonal.  

Condition (ii.c.2): (Verlinde identity) From the definition of the fusion 
algebra associated with the Turaev-Viro-Ocneanu TQFT, $N_{ij}^k$ is a structure 
constant of the fusion algebra. 

From the unitarity of $S$, it is enough to prove 
\begin{equation}\label{verlinde-id}
 \sum_{i,j} N_{ij}^k \; \overline{S_{im}}\; \overline{S_{jn}} = \delta_{mn} 
 \frac{\overline{S_{mk}}}{S_{m0}}.
\end{equation}
The first term of the left-hand side of the equation (\ref{verlinde-id}) is 
written as $Z_\Delta (D^2 \times S^1 \setminus ({\rm Int}D_1^2 
\times S^1 \cup {\rm Int}D_2^2 \times S^1); p_i,p_j;p_k)$. The second and 
the third terms of the left-hand side of the equation (\ref{verlinde-id}) can be 
written as $Z_\Delta (D^2 \times S^1 \setminus {\rm Int }D_0^2 \times S^1;
p_i, S(p_m))$, $Z_\Delta (D^2 \times S^1 \setminus {\rm Int}D_0^2 \times S^1;
p_j, S(p_n))$, respectively. 
Use gluing and embed the second and third terms into the first term, then 
the left-hand side of the equation (\ref{verlinde-id}) is equal to 
$Z_\Delta (D^2 \times S^1 \setminus ({\rm Int}D_1^2 \times S^1 \cup 
{\rm Int}D_2^2 \times S^1); S(p_m),S(p_n);p_k)$. Then, cut $D^2 \times S^1 
\setminus ({\rm Int}D_1^2 \times S^1 \cup {\rm Int}D_2^2 \times S^1$ along 
the meridian. 

\begin{figure}[htbp]
  \begin{center}
\scalebox{0.45}{
\includegraphics{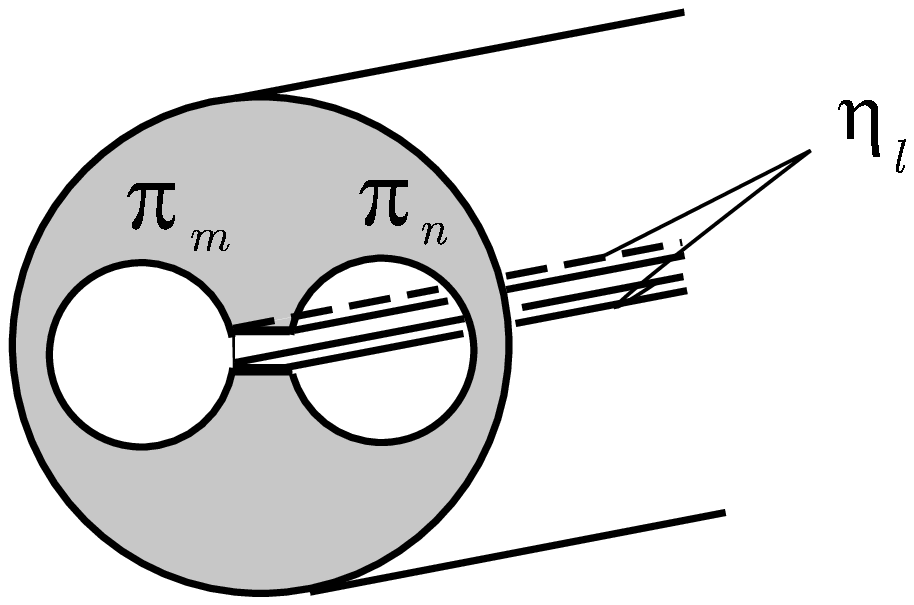}}  
 \end{center}
 \caption{}
 \label{pi_m-pi_n}
\end{figure}

By using the gluing axiom of TQFT, we separate the rectangular part of 
this cut object by taking the summation of $Z_\Delta$ of the 3-manifold 
depicted in Figure \ref{pi_m-pi_n} over orthonormal basis $\eta_l$'s 
on the rectangle. Then, by the identity (\ref{p_i-pi_i}), we get the value 0 
if $m \ne n$ because $\pi_m$ and $\pi_n$ are central orthogonal projections. 
Hence, 
\begin{eqnarray*}
&&Z_\Delta (D^2 \times S^1 \setminus ({\rm Int}D_1^2 \times S^1 \cup 
{\rm Int}D_2^2 \times S^1); S(p_m),S(p_n);p_k) \\
&=& \delta_{mn} \frac{\sum_m Z_\Delta (D^2 \times S^1 \setminus ({\rm Int}D_1^2 
\times S^1 \cup{\rm Int} D_2^2 \times S^1); S(p_m),S(p_n);p_k) \; 
Z_\Delta (D^2 \times S^1;S(p_m))}{Z_\Delta (D^2 \times S^1;S(p_n))} \\
&=& \delta_{mn} \frac{Z_\Delta (D^2 \times S^1 \setminus {\rm Int}D_0^2 \times 
S^1;S(p_m), p_k)}{Z_\Delta (D^2 \times S^1;S(p_n))} \\
&=& \delta_{mn} \frac{\overline{S_{mk}}}{S_{n0}}. 
\end{eqnarray*}

Condition (i): From the Verlinde identity and the unitarity of $S$, we have 
\[
 N_{0j}^k=\sum_{l=0}^r \frac{S_{0l} S_{jl} \overline{S_{kl}}}{S_{0l}} 
 = \sum_{l=0}^r S_{jl} \overline{S_{kl}} = \delta_{jk}.
\]
So $p_0$ is the identitiy element in the fusion algebra. 

Condition (ii.b): Since $S^2$ is $*$-antiisomorphism, it is clear that 
$S^2(p_i)$'s are minimal projections again. 

We have $S^2(p_0)=p_0$ because the definition of $p_0$ is invariant 
under the 180 degree rotation. 

Condition (iv): We note that $Z(U)$ is nothing but the $*$-operation of 
${\rm Tube}\; \Delta$. So $Z(U)(p_i)=p_i^*=p_i$. The canonical map $\theta$ maps 
an orthonormal basis in $Z(-S^1 \times S^1)$ to the dual basis in $Z(S^1 
\times S^1)^*$. So, $\theta(p_i)=\widehat{p_i}$, where 
$\{ \widehat{p_i} \}_{i=0}^r$ is a dual basis of $\{ p_i \}_{i=0}^r$ in 
$Z(S^1 \times S^1)^*$ such that $\widehat{p_i}(p_j)=\delta_{ij}$. Getting 
together, we have $\theta \circ Z(U) (p_i) = \theta(p_i)=\widehat{p_i}$. 
\hfill Q.E.D.

\section{Applications}

From the conclusion of Section 4, we know that there exists a Verlinde basis 
of $Z_{\Delta}(S^1 \times S^1)$ in the sense of Section 2 in a  
Turaev-Viro-Ocneanu TQFT. Hence, for a closed 3-manifold $M$, we have 
the following Dehn surgery formula: 
\[
 Z_\Delta (M) = \sum_{i_1, \dots, i_m=0}^r S_{0i_1} \cdots S_{0i_m} 
 J(L;i_1, \dots, i_m),
\]
where we have assumed that the manifold $M$ is obtained from $S^3$ by Dehn 
surgery along a framed link $L=L_1 \cup \cdots \cup L_m$. 

The purpose in this section is to understand the formula of 
right-hand side of the above equation by introducing a notion of 
the {\it tube system} due to Ocneanu \cite{Ocneanu,O}, which gives a tensor 
category. 

\subsection{Tube systems}

\begin{figure}[htbp]
  \begin{center}
\scalebox{0.47}{
\includegraphics{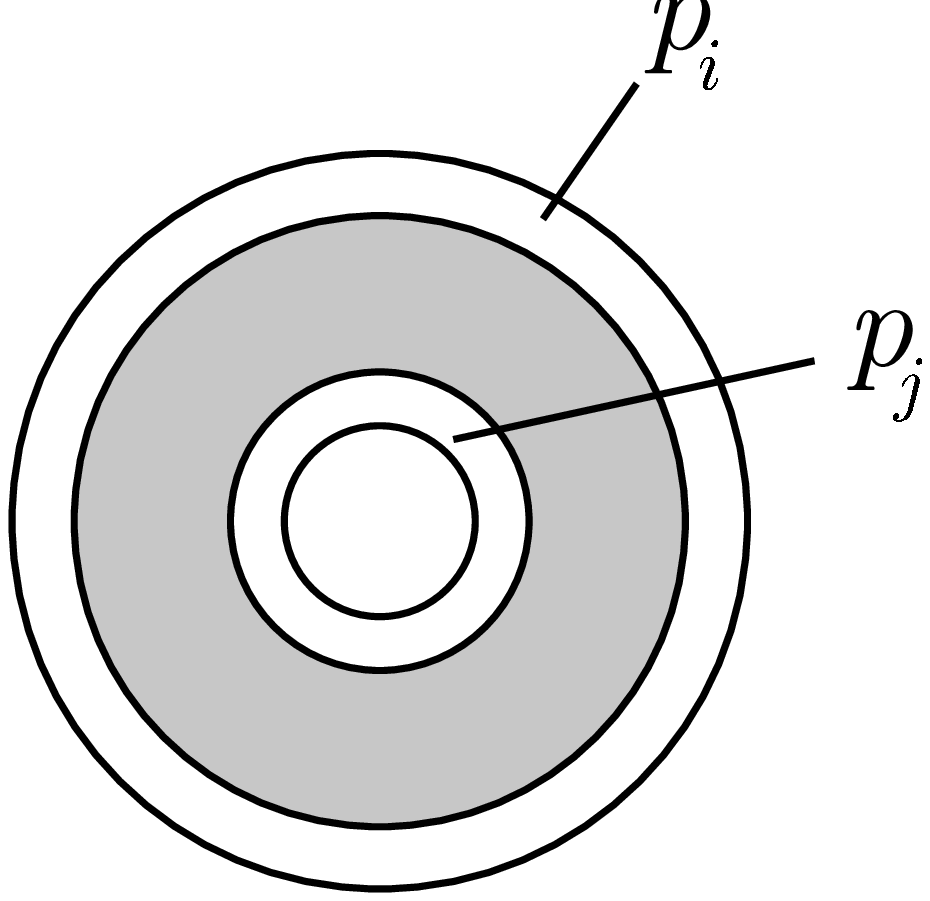}}  
 \end{center}
 \caption{Hom$(p_i,p_j)$}
 \label{p_i-p_j}
\end{figure}

Let us start with the definition of a tube system. Let $\Delta$ be a 
finite system. A {\it tube system} ${\cal D}(\Delta)$ is a tensor category 
defined in the following way. 
First of all, the objects of ${\cal D}(\Delta)$ are the ${\mathbb C}$-linear 
span of all minimal projections in the tube algebra Tube\/$\Delta$. 
For minimal projections $p_i$ and $p_j$, the hom-set ${\rm Hom}(p_i, p_j)$ 
is the set of vectors in the labeled surface depicted in  Figure \ref{p_i-p_j}. 
For general objects $X, Y \in {\cal D}(\Delta)$, we define ${\rm Hom} (X,Y)$ by 
extending ${\rm Hom}(p_i, p_j)$ by linearity. 
Let  $\sim$ be the 
Murray-von Neumann equivalence relation between projections. (Namely, 
two projections $p \; ,q \in {\rm Tube}\; \Delta$ are equivalent in the sense 
of Murray-von Neumann if there exists a partial isometry $v \in {\rm Tube} \; 
\Delta$ such that $p=v^* v$ and $q=v v^*$.) If $p_i \sim p_j$, then 
${\rm Hom}(p_i, p_j) \cong \mathbb{C}$ and if $p_i \nsim p_j$, then 
${\rm Hom}(p_i, p_j) =\{ 0 \}$.  Hence, minimal projections are  
simple objects in ${\cal D}(\Delta)$. 
We denote ${\rm Hom}(p_i, p_i)$ by ${\rm End}(p_i)$ for simplicity. 

\begin{figure}[htbp]
  \begin{center}
\scalebox{0.4}{
\includegraphics{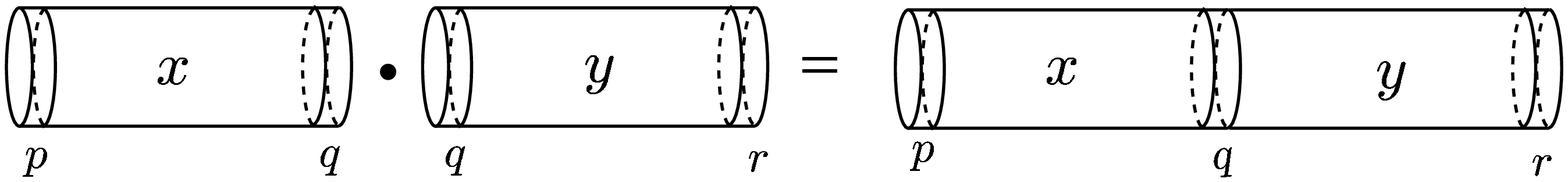}}  
 \end{center}
 \caption{The composition of $x \in {\rm Hom} (q,r)$ and $y \in {\rm Hom}(p,q)$}
 \label{composition}
\end{figure}

For simple objects $p,q$ and $r$, the composition $x \cdot y$ of $x \in 
{\rm Hom} (p,q)$ and $y \in {\rm Hom}(q,r)$ is defined by the concateneation 
of $x$ followed by $y$. See Figure \ref{composition}. For general objects $p,q$ 
and $r$, we define the composition of two morphisms as above by linearlity. 
We also define $x^*$ for $x$ in ${\rm Hom}(p,q)$ by inversing $x$ inside 
out. 

Let $p_j$ be a simple object and $q$ be an object in ${\cal D}(\Delta)$. We 
define an inner product $\langle x, y \rangle$ of $x, y \in {\rm Hom}(q, p_i)$ by 
$\langle x, y \rangle= x \cdot y^* $, as a composition of morphisms. Then, with this 
inner product, ${\rm Hom}(q, p_i)$ becomes a Hilbert space. 

\begin{figure}[htbp]
  \begin{center}
\scalebox{0.65}{
\includegraphics{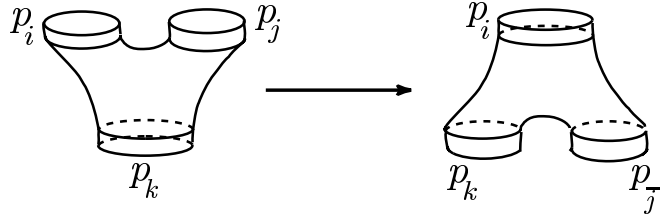}}  
 \end{center}
 \caption{A Frobenius reciprocity map}
 \label{frobenius}
\end{figure}

Next, we define the tensor product $p_i \otimes p_j$ of two simple objects 
$p_i$ and $p_j$ by the fusion product $p_i * p_j$ which was defined in 
Section 2. 
Then, ${\rm Hom} (p_i \otimes p_j, p_k)$ consists of  vectors in the labeled 
surface of Figure \ref{Nijk}. We define an inner product $\langle x, y 
\rangle$ of $x, y \in {\rm Hom} (p_i \otimes p_j, p_k)$ by a composition of two 
morphisms $\langle x, y \rangle= x \cdot y^*  \in {\rm End}(p_k) \cong \mathbb{C}$. 
(Here, $y^*$ is defined by inversing $y$ (pants) inside out.) Note that we have 
Frobenius reciprocities for morphisms of ${\rm Hom}(p_i \otimes p_j, p_k)$. 
For instance, we obtain ${\rm Hom}(p_i \otimes p_j, p_k) \cong {\rm Hom}(
p_i, p_k \otimes p_{\bar{j}})$. See Figure \ref{frobenius}. Frobenius 
reciprocities are given by the graphical operations in the category 
${\cal D}(\Delta)$. 

We call the above defined semi-simple tensor category ${\cal D}(\Delta)$ 
a {\it tube system}. 

\begin{figure}[htbp]
  \begin{center}
\scalebox{0.6}{
\includegraphics{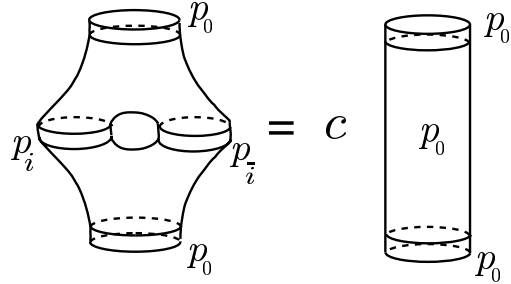}}  
 \end{center}
 \caption{Computing a quantum dimension}
 \label{qdim}
\end{figure}

To compute the quantum dimension of $p_i$, we make the composition of 
morphisms $b_i^* \circ b_i$, where $b_i$ is the distinguished morphism 
from $p_0$  to $p_i \otimes p_{\bar{i}}$ obtained from Frobenius 
reciprocity of $id_{p_i} \in {\rm End}(p_i)$. It must be a scalar multiple 
of $p_0$ since $p_0$ is a simple object. We put this value $c$. 
(See Figure \ref{qdim}.) 

To obtain the value $c$, we connect the upper $p_0$ and the lower $p_0$, 
and embed it into $S^3$. Fill out the outside of the tube, and take the 
Turaev-Viro-Ocneanu invariant of it. Then, we get the 
$c=J(\text{\LARGE $\bigcirc $}; i)/S_{00}$. This is the quantum dimension of 
$p_i$, denoted by ${\rm dim}\; p_i$.

\begin{figure}[htbp]
  \begin{center}
\scalebox{0.8}{
\includegraphics{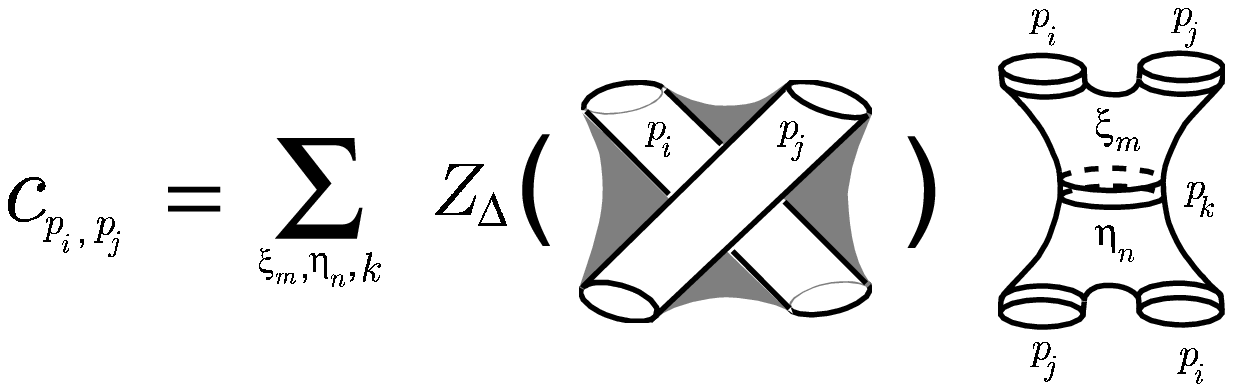}}  
 \end{center}
 \caption{A map $c_{p_i,p_j}$}
 \label{braiding}
\end{figure}

We define the map $c_{p_i,p_j}$ from $p_i \otimes p_j$ to 
$p_j \otimes p_i$ as in Figure \ref{braiding}.
It is easy to see that this map $c_{p_i,p_j}$ satisfies the axioms of 
a braiding, since the map $c_{p_i,p_j}$ is defined in a topological 
way. So, we now know the tensor category ${\cal D}(\Delta)$ is braided. 

\begin{figure}[htbp]
  \begin{center}
\scalebox{0.8}{
\includegraphics{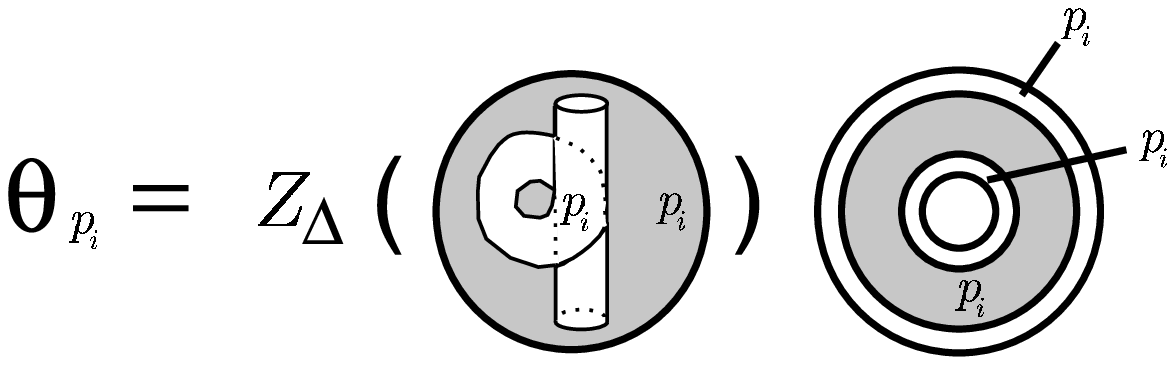}}  
 \end{center}
 \caption{The twist map}
 \label{twist}
\end{figure}

We further define a map $\theta_{p_i} \in {\rm End}(p_i)$ by Figure 
\ref{twist}, where $Z_\Delta$ is evaluated at the 3-ball removed a twisted 
solid tube and the two tubes in the boudaries of it are both labeled by $p_i$. 
This map satisfies the following two equalities. (By the definition of 
$\theta_{p_i}$, these equalities are proven by making pictures corresponding 
to the formulas in both sides and  using some topological moves.)
\begin{eqnarray*}
& \theta_{p_i \otimes p_j} = c_{p_j,p_i} \circ c_{p_i,p_j} \circ (\theta_{p_i} 
 \otimes \theta_{p_j}), \\
& (\theta_{p_i} \otimes id_{p_{\bar{i}}}) \circ b_i = (id_{p_i} \otimes 
 \theta_{p_{\bar{i}}} ) \circ b_i.
\end{eqnarray*}
Namely, $\theta_{p_i}$ defines a twist on ${\cal D}(\Delta)$ and this makes 
${\cal D}(\Delta)$ a ribbon category.  

\begin{figure}[htbp]
  \begin{center}
\scalebox{0.4}{
\includegraphics{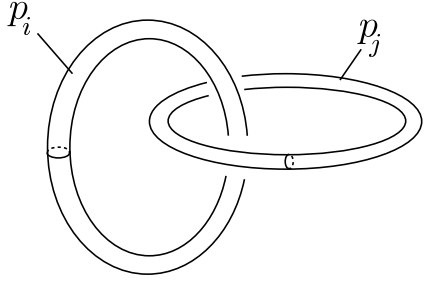}}  
 \end{center}
 \caption{The Hopf link}
 \label{hopf}
\end{figure}

Let us make the following compositions of morphisms in ${\cal D}(\Delta)$. 
\begin{eqnarray*}
&b_0& 
 \circ (id_{p_0} \otimes  b_j) 
 \circ (b_i \otimes id_{p_j} \otimes id_{b_{\bar{j}}}) 
 \circ (id_{p_i} \otimes c_{p_{\bar{i}},p_{j}} \otimes id_{p_{\bar{j}}})
 \circ (id_{p_i} \otimes c_{p_j, p_{\bar{i}}} \otimes id_{p_{\bar{j}}} ) \\
 &\circ& (b_i^* \otimes id_{p_j} \otimes p_{\bar{j}}) 
 \circ (id_{p_0} \otimes b_j^*) 
 \circ b_0^*
\end{eqnarray*}
Then, it is a scalar multiple of $p_0$ and makes a Hopf link $H$ as a 
diagram. (See Figure \ref{hopf}.) 
Let us denote this scalar of the Hopf link $H$ by $s_{ij}$. Embed this 
compositions of tubes into $S^3$ and 
fill out the outside of tubes. By Proposition \ref{Proposition2} in Section 2, 
this provides 
\[
 s_{ij}=\frac{J(H; i, j)}{S_{00}}=\frac{S_{ij}}{S_{00}}, 
\]
where $S=(S_{ij})_{i,j=0}^r$ is the $S$-matrix  with respect to our Verlinde basis 
$\{ p_0, \dots, p_r \}$. 
Since the $S$-matrix is unitary, the matrix $(s_{ij})_{i,j=0}^r$ is 
invertible. It means that our category ${\cal D}(\Delta)$ is modular. 

\begin{rem}
 The notions of the tube algebra and the tube system are also described by 
the language of a II$_1$-subfactor $N \subset M$ with finite Jones index 
and finite depth, although our exposition here uses  an infinite 
subfactor. 
\end{rem}

\subsection{Dehn surgery formula as a Reshetikhin-Turaev invariant}

Let $\cal C$ be a modular category and $\{ V_i \}_{i=0}^r$ its simple 
objects. Put $\Delta=\Delta_+ = \sum_{i=0}^r t_i^{-1} ({\rm dim}\; V_i)^2$, 
$\Delta_- = \sum_{i=0}^r t_i ({\rm dim}\; V_i)^2$ and 
$D=(\sum_{i=0}^r ({\rm dim} V_i)^2)^{1/2}$. (Here, for $\Delta$, we 
followed the notation in \cite{Turaev}.)
Let $M$ be a closed 3-manifold obtained from $S^3$ by Dehn surgery along a 
framed link $L=L_1 \cup \cdots \cup L_m$. Then, the Reshetikhin-Turaev invariant 
of $M$ is given by the formula 
\[
 \tau(M) = \Delta^{\sigma(L)} D^{-\sigma(L)-m-1} \sum_{\lambda \in Col(L)} 
 (\prod_{n=1}^m {\rm dim} V_{\lambda(n)}) F(L,\lambda), 
\]
where $\sigma(L)$ is the signature of $L$ and $F(L,\lambda)$ is 
the invariant of the colored framed link  $(L, \lambda)$. 
(See \cite{Turaev} for the details. Also see \cite{RT}.) 

We note that in the above formula ${\rm dim} V_{\lambda(n)} = 
s_{0\lambda(n)}$. Hence the original Rehshtikhin-Turaev formula can be 
rewritten with the $s$-matrix in the form 
\[
 \tau(M) = \Delta^{\sigma(L)} D^{-\sigma(L)-m-1} \sum_{\lambda \in Col(L)} 
 (\prod_{n=1}^m s_{0 \lambda(n)})  F(L,\lambda).
\]

Now, we start with the modular category ${\cal D}(\Delta)$ defined in 
Section 5.1 instead of a 
general modular category and make the 
Rehsetikhin-Turaev formula. In our case, we already have $S$-matrix 
$( S_{ij} )_{i,j=0}^r$ from the Turaev-Viro-Ocneanu TQFT, which is expressed 
with respect to a Verlinde basis $\{ p_0, \dots, p_r \}$ in 
$Z_\Delta(S^1 \times S^1)$, and since we know that 
$\dim \; p_i=S_{0i}/S_{00}$, by using $D=1/S_{00}$, we can rewrite the 
Reshetikhin-Turaev formula in the form
\[
 \tau(M) = \Delta^{\sigma(L)} D^{-\sigma(L)} \sum_{\lambda \in Col(L)} 
 (\prod_{n=1}^m S_{0 \lambda(n)}) S_{00} F(L,\lambda).
\]

We will prove that $D=\Delta$.
First of all, from the equality $Z_\Delta (S^3)=Z_\Delta (L(1,1))
=\sum_{i=0}^r t_i^{-1} S_{0i}^2$,  we have $\frac{1}{S_{00}}=
\lambda=\sum_{i=0}^r t_i^{-1} S_{0i}^2 $ in our notation. From this, we get 
\[
 \frac{1}{S_{00}}=\lambda =\sum_{i=0}^r t_i^{-1} (\frac{S_{0i}}{S_{00}})^2
 = \sum_{i=0}^r t_i^{-1}({\rm dim} \; p_i)^2=\Delta_+. 
\]
Taking the complex conjuegation of the above formula, we get 
\[
 \lambda = \sum_{i=0}^r t_i ({\rm dim} \; p_i)^2=\Delta_-. 
\]
Thus $\lambda^2 =\Delta_+ \Delta_- = \sum_{i=0}^r 
({\rm dim} \; p_i)^2$. Namely, we have $D=\Delta$. So, the Reshetikhin-Turaev 
invariant constructed from ${\cal D}(\Delta)$ is given by 
\[
 \tau(M) = \sum_{\lambda \in Col(L)}  (\prod_{n=1}^m S_{0 \lambda(n)}) 
 S_{00} F(L,\lambda). 
\]

Since $F$ is uniquely determined by the category of the ribbon tangles 
\cite[Part I, Chapter I]{Turaev}, taking the normalization into consideration, 
we have $F(L,\lambda)=\frac{J(L,\lambda)}{S_{00}}$. Namely, in our 
case, 
\[
 \tau(M) = \sum_{\lambda \in Col(L)}  (\prod_{n=1}^m S_{0 \lambda(n)}) 
 J(L,\lambda), 
\]
which is nothing but our Dehn surgery formula Proposition \ref{Proposition1}. Let us 
summerize this argument as a theorem. 

\begin{thm}\label{main}
 Let $\Delta$ be a finite system. 
 For a closed oriented 3-manifold $M$,  the Reshetikhin-Turaev invariant $\tau(M)$ 
 constructed  from a tube system ${\cal D} (\Delta)$ coincides with the 
 Turaev-Viro-Ocneanu invariant $Z_\Delta (M)$. 
\end{thm}

The following corollary has been proven by several authors \cite{Turaev}, 
\cite{Roberts} etc in various settings. 

\begin{cor}\label{split}
Let $\Delta$ be a finite $C^*$-tensor category arising from a subfactor. 
If  $\Delta$ is a modular category, then ${\cal D} (\Delta)$ is equivalent 
to $\Delta \otimes \Delta^{\rm op}$ and this provides us with $\tau(M) \cdot 
\overline{\tau(M)}=Z_{\Delta}(M)$ for a closed oriented 3-manifold $M$, 
where $\tau$ is the Reshetikhin-Turaev invariant constructed from $\Delta$. 
\end{cor}

\noindent
Proof. In \cite{EK2}, it is proved that ${\cal D}(\Delta)$ is equivalent  
to $\Delta \otimes \Delta^{\rm op}$ when $\Delta$ is a modular category. 
(See Appendix too.) Hence, the rest is clear from Theorem \ref{main}. 
\hfill Q.E.D.

\appendix

\section{Appendix}
In this Appendix, we fix inaccuracies in \cite{EK} and \cite{EK1}.
In \cite{EK} and \cite{EK1}, the authors have analyzed the structure of 
$M_\infty$-$M_\infty$ bimodules obtained from the {\it asymptotic 
inclusion} $M \vee M^{\rm op} \subset M_\infty$ starting from 
the inclusion of AFD II$_1$ factors $N \subset M$ with finite Jones index 
and finite depth, following Ocneanu. 

Let $\{p_0, \cdots, p_r \}$ be a representative set of the equivalence classes of 
minimal projections of Tube$\Delta$ by the Murray-von Neumann equivalence relation. 

We present the $M_\infty$-$M_\infty$ bimodule $X(p_i)$ by Figure \ref{a1}, 
where the circle at the middle is empty and the annulus around it is 
labeled with the minimal projection $p_i$.
(See Chapter 12 in \cite{EK} for more explanation of this kind of pictures.)
\unitlength 0.5mm
\thicklines
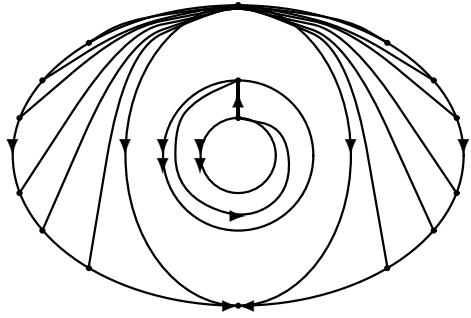
\begin{figure}[htpb]
\begin{center}\begin{picture}(160,100)
\put(70,50){\circle{20}}
\put(70,50){\circle{40}}
\put(70,60){\circle*{1}}
\put(70,70){\circle*{1}}
\spline(70,70)(59,65)(53,58)(53.5,40.5)(71,32)(84,40)(83,57)(70,60)
\put(72,34){\vector(1,0){0}}
\put(70,67){\vector(0,1){0}}
\put(70,60){\line(0,1){10}}
\put(50,50){\vector(0,-1){0}}
\put(50,45){\vector(0,-1){0}}
\put(60,50){\vector(0,-1){0}}
\put(60,45){\vector(0,-1){0}}
\put(70,50){\ellipse{60}{80}}
\put(70,50){\ellipse{120}{80}}
\put(10,50){\vector(0,-1){0}}
\put(40,50){\vector(0,-1){0}}
\put(100,50){\vector(0,-1){0}}
\put(130,50){\vector(0,-1){0}}
\put(70,10){\vector(1,0){0}}
\put(70,10){\vector(-1,0){0}}
\put(70,10){\circle*{1}}
\put(70,90){\circle*{1}}
\put(11.9,40){\circle*{1}}
\put(11.9,60){\circle*{1}}
\put(128.1,40){\circle*{1}}
\put(128.1,60){\circle*{1}}
\put(18,30){\circle*{1}}
\put(18,70){\circle*{1}}
\put(122,30){\circle*{1}}
\put(122,70){\circle*{1}}
\put(30.3,20){\circle*{1}}
\put(30.3,80){\circle*{1}}
\put(109.7,20){\circle*{1}}
\put(109.7,80){\circle*{1}}
\spline(30.3,80)(40,86)(60,89)(70,90)(80,89)(100,86)(109.7,80)
\spline(18,70)(40,84)(60,89)(70,90)(80,89)(100,84)(122,70)
\spline(11.9,60)(40,84)(60,88)(70,90)(80,88)(100,84)(128.1,60)
\spline(11.9,40)(40,83)(60,87)(70,90)(80,87)(100,83)(128.1,40)
\spline(18,30)(40,80)(60,87)(70,90)(80,87)(100,80)(122,30)
\spline(30.3,20)(40,78)(60,87)(70,90)(80,87)(100,78)(109.7,20)
\end{picture}
\end{center}
\caption{The bimodule $X(p_i)$}
\label{a1}
\end{figure}
We can prove that each $M_\infty$-$M_\infty$ bimodule $X(p_i)$ $(i=0, 
\cdots, r)$ is irreducible in a similar way to the proof of 
Theorem 12.26 in \cite{EK}. 

Since $[{}_{M_\infty} X(p_i)_{M_\infty}]^{1/2}=\dim \; p_i$, 
we have the equalities $\sum_{i=0}^r [X(p_i)]=\sum_{i=0}^r (\dim \; p_i)^2
=\lambda^2$, where dim $p_i$ is the quantum dimension of $p_i$ as an object of 
${\cal D}(\Delta)$. On the other hand, the global index of the asymptotic 
inclusion 
is given by $\lambda^2$. This means that all the irreducible 
$M_\infty$-$M_\infty$ bimodules obtained from the asymptotic inclusion are 
given by $\{X(p_i) \}_{i=0}^r$. 

\begin{rem}
In \cite{EK} and \cite{EK1}, irreducible bimodules are labeled by $\pi_i$'s, 
i.e., minimal central projections of Tube$\Delta$, instead of the minimal 
projections of it, which is incorrect.
\end{rem}

We describe the fusion rule of irreducble $M_\infty$-$M_\infty$ 
bimodules. The relative tensor product of two $M_\infty$-$M_\infty$ 
bimodules $X(p_i) \otimes_{M_\infty} X(p_j)$ is decomposed into irreducible 
bimodules as in Figure \ref{a5}.

\unitlength 0.45mm
\thicklines
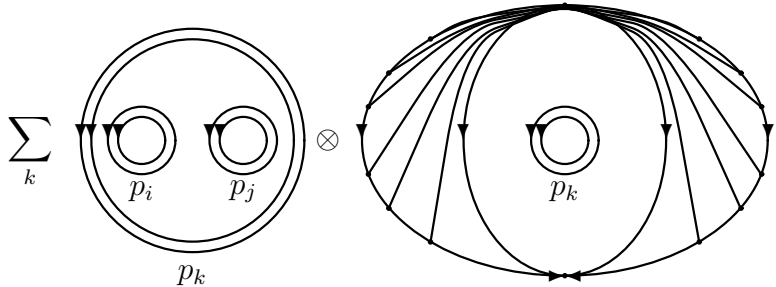
\begin{figure}[htpb]
\begin{center}\begin{picture}(240,100)
\put(170,50){\circle{20}}
\put(170,50){\circle{14}}
\put(45,50){\circle{20}}
\put(45,50){\circle{14}}
\put(75,50){\circle{20}}
\put(75,50){\circle{14}}
\put(60,50){\circle{60}}
\put(60,50){\circle{66}}
\put(160,50){\vector(0,-1){0}}
\put(163,50){\vector(0,-1){0}}
\put(30,50){\vector(0,-1){0}}
\put(27,50){\vector(0,-1){0}}
\put(35,50){\vector(0,-1){0}}
\put(38,50){\vector(0,-1){0}}
\put(65,50){\vector(0,-1){0}}
\put(68,50){\vector(0,-1){0}}
\put(170,50){\ellipse{60}{80}}
\put(170,50){\ellipse{120}{80}}
\put(110,50){\vector(0,-1){0}}
\put(140,50){\vector(0,-1){0}}
\put(200,50){\vector(0,-1){0}}
\put(230,50){\vector(0,-1){0}}
\put(170,10){\vector(1,0){0}}
\put(170,10){\vector(-1,0){0}}
\put(170,10){\circle*{1}}
\put(170,90){\circle*{1}}
\put(111.9,40){\circle*{1}}
\put(111.9,60){\circle*{1}}
\put(228.1,40){\circle*{1}}
\put(228.1,60){\circle*{1}}
\put(118,30){\circle*{1}}
\put(118,70){\circle*{1}}
\put(222,30){\circle*{1}}
\put(222,70){\circle*{1}}
\put(130.3,20){\circle*{1}}
\put(130.3,80){\circle*{1}}
\put(209.7,20){\circle*{1}}
\put(209.7,80){\circle*{1}}
\spline(130.3,80)(140,86)(160,89)(170,90)(180,89)(200,86)(209.7,80)
\spline(118,70)(140,84)(160,89)(170,90)(180,89)(200,84)(222,70)
\spline(111.9,60)(140,84)(160,88)(170,90)(180,88)(200,84)(228.1,60)
\spline(111.9,40)(140,83)(160,87)(170,90)(180,87)(200,83)(228.1,40)
\spline(118,30)(140,80)(160,87)(170,90)(180,87)(200,80)(222,30)
\spline(130.3,20)(140,78)(160,87)(170,90)(180,87)(200,78)(209.7,20)
\put(170,35){\makebox(0,0){$p_k$}}
\put(45,35){\makebox(0,0){$p_i$}}
\put(75,35){\makebox(0,0){$p_j$}}
\put(60,10){\makebox(0,0){$p_k$}}
\put(12,47){\makebox(0,0){$\displaystyle\sum_k$}}
\put(100,50){\makebox(0,0){$\otimes$}}
\end{picture}
\end{center}
\caption{Decomposition of a bimodule}
\label{a5}
\end{figure}
Hence, the fusion rule is given by the fusion rule of the fusion algebra 
associated with the Turaev-Viro-Ocneanu TQFT. 

It is easy to see that $M_\infty$-$M_\infty$ bimodules arising from the 
asymptotic inclusion $M \vee M^{\rm op} \subset M_\infty$ give rise to  
a modular category with the same braiding and twist as ones in a tube system. 
We denote this category of $M_\infty$-$M_\infty$ bimodules by ${\cal M}_\infty$. 
To a simple object $X(p_i)$ in ${\cal M}_\infty$, we assign a simple object 
$p_i$ in ${\cal D}(\Delta)$. 
From Figure \ref{a5}, it is easy to see that we have a functor 
$F: \; {\cal M}_\infty \longrightarrow {\cal D}(\Delta)$, when we look at 
the morphisms in ${\cal M}_\infty$. 

Let us now consider the opposite direction to the functor $F$. 
For a given morphism $x \in {\rm Hom} (p_i \otimes p_j, p_k)$ of 
${\cal D}(\Delta)$, we will construct a homomoprhism in ${\rm Hom} (X(p_i) 
\otimes_{M_\infty} X(p_j), X(p_k))$. 

\unitlength 0.5mm
\thicklines
\begin{figure}[htpb]
\begin{center}\begin{picture}(160,100)
\put(85,50){\circle{20}}
\put(55,50){\circle{20}}
\put(85,50){\circle{14}}
\put(55,50){\circle{14}}
\put(45,50){\vector(0,-1){0}}
\put(75,50){\vector(0,-1){0}}
\put(48,50){\vector(0,-1){0}}
\put(78,50){\vector(0,-1){0}}
\put(70,50){\ellipse{60}{80}}
\put(70,50){\ellipse{120}{80}}
\put(10,50){\vector(0,-1){0}}
\put(40,50){\vector(0,-1){0}}
\put(100,50){\vector(0,-1){0}}
\put(130,50){\vector(0,-1){0}}
\put(70,10){\vector(1,0){0}}
\put(70,10){\vector(-1,0){0}}
\put(70,10){\circle*{1}}
\put(70,90){\circle*{1}}
\put(11.9,40){\circle*{1}}
\put(11.9,60){\circle*{1}}
\put(128.1,40){\circle*{1}}
\put(128.1,60){\circle*{1}}
\put(18,30){\circle*{1}}
\put(18,70){\circle*{1}}
\put(122,30){\circle*{1}}
\put(122,70){\circle*{1}}
\put(30.3,20){\circle*{1}}
\put(30.3,80){\circle*{1}}
\put(109.7,20){\circle*{1}}
\put(109.7,80){\circle*{1}}
\spline(30.3,80)(40,86)(60,89)(70,90)(80,89)(100,86)(109.7,80)
\spline(18,70)(40,84)(60,89)(70,90)(80,89)(100,84)(122,70)
\spline(11.9,60)(40,84)(60,88)(70,90)(80,88)(100,84)(128.1,60)
\spline(11.9,40)(40,83)(60,87)(70,90)(80,87)(100,83)(128.1,40)
\spline(18,30)(40,80)(60,87)(70,90)(80,87)(100,80)(122,30)
\spline(30.3,20)(40,78)(60,87)(70,90)(80,87)(100,78)(109.7,20)
\put(55,35){\makebox(0,0){$p_i$}}
\put(85,35){\makebox(0,0){$p_{\bar{i}}$}}
\end{picture}
\end{center}
\caption{The bimodule $X(p_i) \otimes_{M_\infty} X(p_{\bar{i}})$}
\label{a3}
\end{figure}
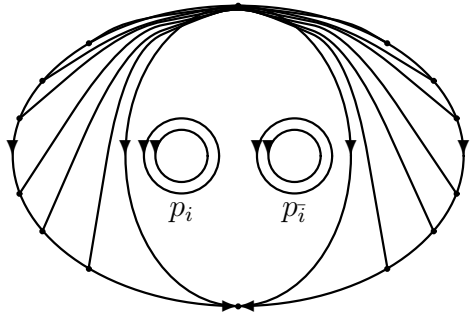

\begin{figure}[htbp]
  \begin{center}
\scalebox{0.45}{
\includegraphics{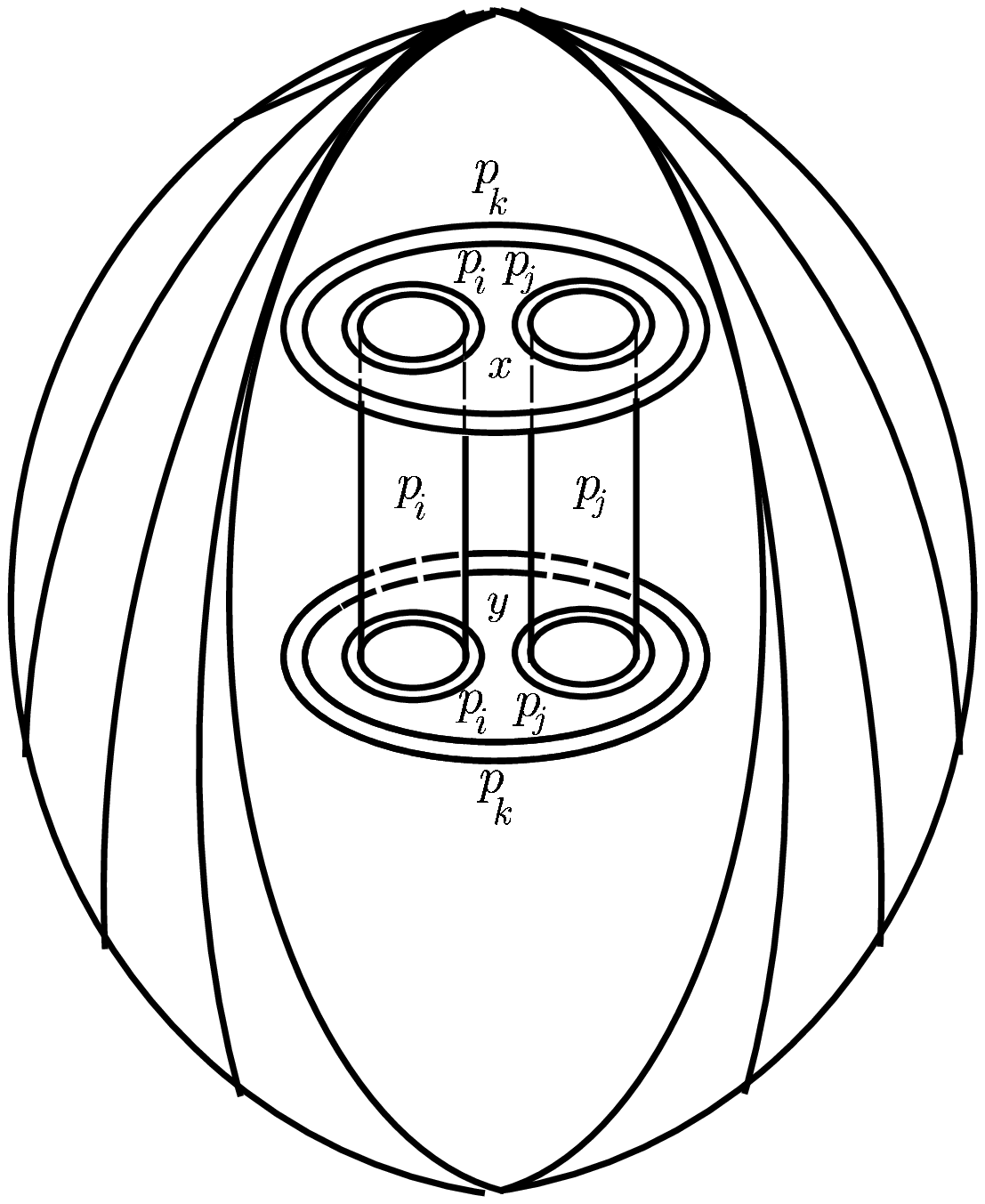}}  
 \end{center}
 \caption{}
 \label{a6}
\end{figure}

Let $y_0$ be an arbitrary element of $X(p_i) \otimes_{M_\infty} X(p_j)$. This 
$y_0$ is considered as an element in the bimodule in Figure \ref{a3}. We 
denote it by $y_1$. Then, we attach $x \in {\rm Hom}(p_i \otimes p_j, p_k)$ to 
this $y_1$ by using two tubes $p_i$, $p_j$. See Figure \ref{a6}. 
\begin{figure}[htbp]
  \begin{center}
\scalebox{0.9}{
\includegraphics{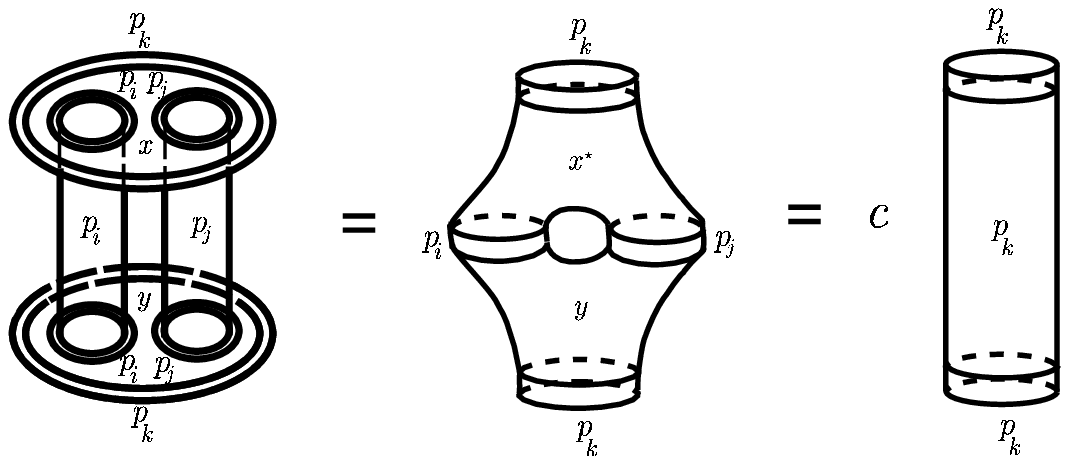}}  
 \end{center}
 \caption{}
 \label{a7}
\end{figure}
The central part of Figure \ref{a6} can be viewed as in Figure \ref{a7}. 
Hence, we get a morphism in ${\rm Hom}(X(p_i) \otimes_{M_\infty} X(p_j), 
X(p_k))$. This induces a functor $G: \; {\cal D}(\Delta) \longrightarrow 
{\cal M}_\infty$. It is not difficult to see that $F$ and $G$ are 
functors, which preserve the operations in modular categories, and inverse to 
each other. 

Hence, the category of $M_\infty$-$M_\infty$ bimodules obtained from the 
asymptotic inclusion ${\cal M}_\infty$ is equivalent to the tube 
system ${\cal D}(\Delta)$ constructed from $N \subset M$ as modular 
categories.

\end{document}